\documentclass[reqno,9pt,article]{amsart}
\usepackage{amsmath,amsxtra,amssymb,latexsym, amscd,amsthm}
\usepackage[unicode]{hyperref}
\usepackage{array,tabularx,longtable,multicol,indentfirst,fancybox,color}
\usepackage{cases}
\usepackage{graphicx}
\usepackage{multicol}
\usepackage{mathrsfs}
\usepackage{enumerate}
\usepackage{tikz}
\usepackage{tikz-cd} 
\usepackage{tikz,tkz-tab}
\usetikzlibrary{calc, angles, quotes, intersections, positioning, patterns}

\advance\hoffset-1truecm\relax
\newtheorem{proposition}{Proposition}[section]

\newtheorem{theorem}{Theorem}[section]
\newtheorem{definition}{Definition}
\newtheorem{lemma}{Lemma}[section]
\newtheorem{corollary}{Corollary}[section]
\newtheorem{remark}{Remark}
\newtheorem{example}{Example}[]
\numberwithin{equation}{section}

\def\R{\mathbb{R}}
\def\C{\mathbb{C}}
\def\a{\alpha}

\begin{document}
	\title[A refined variant of Hartley convolution: Algebraic structures, spectral radius and related issues]{A refined variant of Hartley convolution:\\ Algebraic structures, spectral radius and related issues}
	\author[Trinh Tuan]{Trinh Tuan}
	\date{01 February 2026, 
		Accepted by \textsc{Integral Transforms and Special Functions}}
\thanks{For a Special Issue on
	Classical Analysis and Beyond: 25 Years of Semyon Yakubovich's Academic Excellence in Portugal}	 
\thanks{Published online: 17 February 2026.  \url{https://doi.org/10.1080/10652469.2026.2626960}}

	\maketitle
	\begin{center}
	Department of Mathematics, Faculty of Natural Sciences, Electric Power University,\\ 235-Hoang Quoc Viet Rd., Hanoi, Vietnam.\\		
	E-mail: \texttt{tuantrinhpsac@yahoo.com}
	\end{center}
			
\begin{abstract}
In this work, we propose a novel convolution product associated with the $\mathscr{H}$-transform, denoted by $\underset{\mathscr{H}}{\ast}$,  and explore its fundamental properties. Here, the $\mathscr{H}$-transform may be regarded as a refined variant of the classical Fourier, Hartley transform, with kernel function depending on two parameters $a,b$. Our first contribution shows that the space of integrable functions, equipped with multiplication given by the $\underset{\mathscr{H}}{\ast}$-convolution, constitutes the commutative Banach algebra over the complex field, albeit without an identity element. Second, establishes the Wiener--Lévy type invertibility criterion for $\mathscr{H}$-algebras, obtained through the density property and process of unitarization, which serves as a key step toward the proof of Gelfand's spectral radius theorem. Third, provides an explicit upper-bound of Young's inequality for $\underset{\mathscr{H}}{\ast}$-convolution and its direct corollary. Finally, all of these theoretical findings are applied to analyze specific classes of the Fredholm integral equations and heat source problems, yielding a priori estimates under the established assumptions.

\end{abstract}	

\section{Introduction}\label{sec1}
\subsection{Background}The theory of convolution in integral transforms has long been a vibrant and actively pursued area of research in applied mathematics, engineering, and physics \cite{Poularikas1996,BGM1985}. 
The behavior of a transform with respect to convolution depends both on the specific form of the convolution and on the associated transform. 
A typical
example arises in the case of the Hankel convolution, where the transform of the convolution of two functions is
given by the product of their transforms, multiplied by a power function determined by the order of the Hankel
transform (see Vu Kim Tuan and Megumi Saigo \cite{TuanSaigo95}). In many situations, the existence of a convolution can be inferred from the existence of the corresponding transforms of the underlying functions, even when the transform in question is not directly linked to the convolution formula.
Moreover, the convolution can be extended beyond functions on the complex plane. Indeed, a generalized integral transform was established by Hyun Soo Chung and Vu Kim Tuan in~\cite{ChungTuan11}, for a class of functionals, specifically the space of complex-valued continuous functions on a finite interval $[0,T]$ vanishing at the origin. While convolution operations are often defined in relation to specific transforms, it is noteworthy that such an association is not always necessary. 
However, the concept which is central to integral transforms is still governed by the analytical properties of its kernel function, with the Fourier and Hartley transforms serving as the most well-known examples, together with the convolution structures naturally associated with them and applications~\cite{TuanVK92,TMY86,Yaku94Luchko,tuan22,MMA2022,Tuan25}.
We begin by recalling some notions on the classical Fourier and Hartley integral transforms, which were introduced in \cite{Titchmarsh1986,Bracewell86}. The Fourier transform of the function $f$, denoted by \( (F) \), is defined by
\begin{equation}\label{eq1.1}
(Ff)(y) = (2\pi)^{-n/2} \int_{\mathbb{R}^n} e^{i x y} f(x)\, dx, \ y \in \mathbb{R}^n,
\end{equation}
and its corresponding reverse transform is given by the formula
$
f(x)=(F^{- 1} f)(y) = (2\pi)^{-n/2} \int_{\mathbb{R}^n} e^{-i x y} f(y)\, dy.
$ The two-sided Fourier sine and cosine integral transforms, denoted by \( (F_s) \) and \( (F_c) \), respectively, are defined by integral formulas
\begin{equation}\label{FsFc}
(F_s f)(y) = (2\pi)^{-n/2} \int_{\mathbb{R}^n} \sin(x y) f(x) dx \quad \text{and} \quad 
(F_c f)(y) = (2\pi)^{-n/2} \int_{\mathbb{R}^n} \cos(x y) f(x) dx.
\end{equation}
The Hartley transform is denoted by \( (H) \), which is defined by
\begin{equation}\label{eq1.2}
(Hf)(y)=(2\pi)^{-n/2} \int_{\mathbb{R}^n} \mathrm{cas}(xy) f(x) dx,\ y\in \mathbb{R}^n,
\end{equation} where $\operatorname{cas}(xy)=\cos(xy)\pm\sin(xy)$ as classical Hartley kernel,  $x=(x_1,x_2,\dots,x_n) \in \R^n$, $y=(y_1,y_2,\dots,y_n) \in \R^n$ with $xy=\sum_{i=1}^n x_i y_i$ and $dx=dx_1 dx_2 \dots dx_n$. A recent result in \cite{Castro22} by Castro, Guerra, and Tuan introduced  an alternative formulation for the classical Hartley kernel by considering the following integral transform
\begin{equation}\label{eq1.3}
(\mathscr{H}f)(y):=(2\pi)^{-n/2} \int_{\mathbb{R}^n} (a\cos xy+b\sin xy)f(x) dx,\ y\in\mathbb{R}^n,
\end{equation}
where \(a, b \in \R\) are non-zero constants.  
For the convenience of unifying notations, throughout in this paper we will refer to \eqref{eq1.3} as the $\mathscr{H}$-transform. We now make some observations to clarify that the $\mathscr{H}$-transform can be regarded as a refined variant of the classical Fourier and Hartley transforms. 
Obviously, it can be represented as a linear combination of the Fourier cosine and sine transforms, namely:
$
\mathscr{H} = a (F_c) + b (F_s).
$
Moreover, when \( a = 1 \) and \( b = i \), then  $\mathscr{H}$-transform \eqref{eq1.3} coincides with the classical Fourier transform $( F )$ as defined in~\eqref{eq1.1}. In the special case where \( a = 1 \) and \( b = \pm 1 \), then $\mathscr{H}$-transform \eqref{eq1.3} reduces to the classical Hartley transform \eqref{eq1.2}. Furthermore, due to  Theorem 2.3 in \cite{Castro22} asserts that:  If $f \in L_1(\mathbb{R}^n)$ and $(\mathscr{H}f)(y)\in L_1(\mathbb{R}^n)$, then the reverse transform of \eqref{eq1.3} defined by 
\begin{equation}\label{eq1.4}
(\mathscr{H}^{-1}f)(y):= (2\pi)^{-n/2}\int_{\mathbb{R}^n}\left(\frac{1}{a}\cos xy+\frac{1}{b}\sin xy\right) f(y) dy,\ x\in \mathbb{R}^n.
\end{equation}
Besides, it should also be noticed that, the $\mathscr{H}$-transform admit representations in terms of the classical Fourier transform $(F)$ and its reverse $(F^{-1})$, given explicitly by
$\mathscr{H}=\frac{a-ib}{2} F + \frac{a+ib}{2} F^{-1}.
$
Denote by \(\mathscr{S}(\mathbb{R}^n)\) the Schwartz space \cite{AdamsFournier2003sobolev}, consisting 
of all rapidly decreasing smooth functions on \(\mathbb{R}^n\). 
It is well known that \(\mathscr{S}(\mathbb{R}^n)\) is dense in both 
\(L_1(\mathbb{R}^n)\) and \(L_2(\mathbb{R}^n)\). Since the mapping \(f \mapsto Ff\) is continuous (with respect to the \(L_2\)-metric) of the dense subspace \(\mathscr{S}(\mathbb{R}^n) \subset L_2(\mathbb{R}^n)\) into \(\mathscr{S}(\mathbb{R}^n)\), it admits a unique continuous extension   \(\overline{F}: L_2(\mathbb{R}^n) \to L_2 (\mathbb{R}^n)\).
Due to Theorem~7.9 in \cite{Rudin1991functional}, there exists a unique linear isomorphism \(\overline{F}\) of \(L_2(\mathbb{R}^n)\) onto itself, uniquely determined by the requirement that \((\overline{F} f)(y) = (F f)(y)\) for all \(f \in \mathscr{S}(\mathbb{R}^n)\). Similarly, its reverse operator \(\overline{F^{-1}}\) is uniquely determined by the condition \((\overline{F^{-1}} f)(y) = (F^{-1} f)(y)\) for every \(f \in \mathscr{S}(\mathbb{R}^n)\). In particular, $
(\overline{F}f)(y) = \lim\limits_{R\to \infty} (2\pi)^{-n/2}
\int_{|x|\leq R} e^{-i x y}\, f(x)\, dx.
$  Furthermore, Theorem 2.2 in \cite{Castro22} has proven that the $\mathscr{H}$-transform is a linear operator from \(\mathscr{S}(\mathbb{R}^n)\) into itself, possesses the reflection property, and is invertible on \(\mathscr{S}(\mathbb{R}^n)\). Therefore, for any functions \(f \in L_2(\mathbb{R}^n)\), then both \(\mathscr{H}f\) and \(\mathscr{H}^{-1}f\) belong to \(L_2(\mathbb{R}^n)\). However
, recall that both the classical Fourier and Hartley transforms are unitary operators on 
\(L_2(\mathbb{R}^n)\), i.e., they preserve the inner product and the 
$L_2$-norm   
for all $ f \in L_2(\mathbb{R}^n).
$  
In contrast, the $\mathscr{H}$-transform does not, in general, preserve the 
$L_2$-norm, since the coefficients \(a, b\) in its kernel alter the 
normalization. This lack of norm-preservation is precisely the reason why 
$\mathscr{H}$-transform  is termed a non-unitary operator on the 
Hilbert space \(L_2(\mathbb{R}^n)\).  This also constitutes the most evident 
difference when compared with the  classical Fourier and Hartley transforms. \\ 
\noindent  Now we recall the concept of a complex Banach algebra. 
Let $X$ be a complex vector space equipped with a norm $\|.\|_X$. 
Suppose that $X$ is complete with respect to this norm, so that $(X,\|.\|_X)$ is a Banach space. 
If, in addition, $X$ is equipped with a multiplication ``$\cdot$'' satisfying:  

\begin{itemize}
	\item associativity: $(xy)z = x(yz)$,
	\item distributivity: $(x+y)z = xz + yz$ and $x(y+z) = xy + xz$,
	\item scalar compatibility: $a(xy) = (a x)y = x(a y)$, $\forall x,y,z \in X, \ a \in \mathbb{C}$,
\end{itemize}
\noindent together with the submultiplicative property:
$
\|xy\|_X \leq \|x\|_X \, \|y\|_X, \forall x,y \in X,
$
then $(X,+,\cdot,\|.\|_X)$ is called a \textit{Banach algebras} \cite{Dales00,Kaniuth09}.  
More generally, if $X$ is a normed vector space over $\mathbb{C}$ endowed with an associative multiplication satisfying the submultiplicative inequality above (but not necessarily complete), then $X$ is called a \textit{Normed algebra} (sometimes referred to as a \emph{Normed ring}, a terminology originating from Gelfand's thesis).  One of Gelfand's greatest achievements was the theory of commutative normed rings, which he created and investigated in his Doctor of Science thesis submitted in 1938. Together with his subsequent work \cite{Gelfand-Naimark} on non-commutative rings of linear operators on Hilbert spaces in 1943, these studies constituted decisive steps in the development of functional analysis and exerted a profound influence on related fields such as algebraic geometry and theoretical physics. The significance of these works is emphasized in \cite{50namGelfand}: it was Gelfand who brought to light the fundamental concept of a maximal ideal, thereby enabling the unification of previously disconnected results and laying the foundation for a new and profound theory. Moreover, Gelfand’s theory of normed rings revealed close connections between Banach's general functional analysis and classical analysis \cite{1964,Naimark}.
In the present work, we focus on the space 
$
X \equiv L_1(\mathbb{R}^n),
$
the space of Lebesgue integrable functions on $\mathbb{R}^n$. 
The setting up of a Banach algebra over the complex field may then be viewed schematically as:  
$
\text{Banach space}\rightarrow\text{Normed algebra (Normed ring)}\rightarrow \text{Banach algebra}.
$ 
It is well known that $L_1(\mathbb{R}^n)$, equipped with the $L_1$-norm, forms a Banach space. 
Therefore, a natural question that arises here: 
Does there exist a convolution product defined via the $\mathscr{H}$-transform such that, when regarded as the multiplication on $L_1(\mathbb{R}^n)$, it turns this space into a Banach algebra structure?
Addressing this question constitutes the central motivation of our work. As a first step toward this goal, we propose the following definition of the convolutional structure associated with the $\mathscr{H}$-transform:

\begin{definition}[Main object]\label{def2.1}
	The convolution of two functions  $f$, $g$ associated with the $\mathscr{H}$-transform \eqref{eq1.3} is denoted by $(f \underset{\mathscr{H}}{\ast} g)$ and defined by the formula
	\begin{equation}\label{eq2.1}
	\big(f \underset{\mathscr{H}}{\ast} g\big)(x) :=\frac{1}{4a(2\pi)^{\frac{n}{2}}} \int_{\mathbb{R}^n} K_{(a,b)}[f](x,v) g(v) dv,\ \ x\in\mathbb{R}^n,
	\end{equation}
	where the kernel $K_{(a,b)}[f](x,v)$ is represented by
	\begin{equation}\label{eq2.2}
	K_{(a,b)}[f](x,v)=(3a^2-b^2)f(x-v) + (a^2+b^2)\big[f(x+v)+f(-x+v)-f(-x-v)\big].
	\end{equation}
\end{definition}
\subsection{Organization} This paper is organized into five sections. In Section \ref{botro}, we prove some auxiliary lemmas. First, we discuss the validity of the conditions to ensure the non-degeneracy of kernel function in \eqref{eq1.3} and demonstrate a lemma for the injectivity of $\mathscr{H}$-transform. The second lemma establishes a Riemann--Lebesgue type lemma for \eqref{eq1.3}, while the third lemma provides an estimate for the kernel \eqref{eq2.2} associated with $\underset{\mathscr{H}}{\ast}$-convolution. All three play an essential role in preparing the ground for the proofs of the main results later on. As outlined in the Abstract, the primary contributions of this paper are presented in Section~\ref{sec2} and Section~\ref{sec3}. Finally,  Section~\ref{sec4} is devoted to exploring the applicability of the convolution structure~\eqref{eq2.1} to the solvability of certain heat source problems and Fredholm integral equations. This analysis relies on the theoretical framework developed in the preceding sections and makes essential use of the Wiener–Lévy invertibility criterion for $\mathscr H$-algebras, and together with related convolution inequalities.
\section{Some auxiliary lemmas}\label{botro}
\noindent First, we will analyze the validity of the conditions for injectivity of the $\mathscr{H}$-transform \eqref{eq1.3}. $\mathscr{H}$ admits the representation $
\mathscr{H}f = \frac{a - ib}{2} F f + \frac{a + ib}{2} F^{-1} f,
$
where $F$ denotes the classical Fourier transform with normalization $(2\pi)^{-n/2}$ symmetric and $F^{-1}f(y) = Ff(-y)$. Now, assume $\mathscr{H}f = 0$ a.e. for some $f \in L_1(\mathbb{R}^n)$. Setting $g = Ff$, we obtain the functional equation
\begin{equation}\label{**}
(a - ib) g(y) + (a + ib) g(-y) = 0 \quad \text{a.e. in }\  \mathbb{R}^n.
\end{equation}
If $a + ib \neq 0$, then equation \eqref{**} yields
$
g(-y) = -\frac{a - ib}{a + ib} \, g(y).
$
Substituting $-y$ for $y$, we have
$$
g(y) = -\frac{a - ib}{a + ib} \, g(-y) = \left(-\frac{a - ib}{a + ib}\right)^2 g(y).
$$
Denote $k = -\frac{a - ib}{a + ib}$. Then $g(y) = k^2 g(y)$ equivalent to
$
(1 - k^2) g(y) = 0$ for a.e. $y$, and we infer that
\[
k^2 = \frac{(a - ib)^2}{(a + ib)^2} = 1 \Leftrightarrow (a - ib)^2 = (a + ib)^2 \Leftrightarrow 4iab = 0.
\]
Thus $k^2 = 1$ occurs precisely when $ab = 0$, i.e., when $a = 0$ or $b = 0$. Therefore: 
\begin{itemize}
	\item If $a \neq 0$ and $b \neq 0$, then $k^2 \neq 1$, and equation \eqref{**} forces $g = 0$ a.e. 
	Since the Fourier transform $F$ is injective on $L_1(\mathbb{R}^n)$, we conclude $f = 0$ almost everywhere. 
	Hence $\mathscr{H}$ is injective.
	
	\item If $a = 0$ and $b \neq 0$, then $\mathscr{H}f = b F_s f$, where $F_s$ is the two-sided Fourier sine transform. 
	Since $F_s$ annihilates all even functions, this implies $\mathscr{H}$ is not injective.
	
	\item If $a \neq 0$ and $b = 0$, then $\mathscr{H}f = a F_c f$, where $F_c$ is the two-sided Fourier cosine transform. 
	Since $F_c$ annihilates all odd functions, this implies $\mathscr{H}$ is not injective.
\end{itemize}
\noindent On the other hand, it's not difficult to recognize that the condition $a \neq \pm b$  is not necessary for injectivity per se. For instance, when $a = b \neq 0$, we have $k = -\frac{1-i}{1+i} = -i$, so $k^2 = -1 \neq 1$, and injectivity still holds. 
Similarly for the case $a = -b \neq 0$. Thus, we have proved
\begin{lemma}[The injectivity of $\mathscr H$-transform]\label{donanh}
	Let $a, b \in \mathbb{R}$ with $a \neq 0$ and $b \neq 0$. 
	Then  $\mathscr{H} : L_1(\mathbb{R}^n) \to C_0(\mathbb{R}^n)$ defined by \eqref{eq1.3} is injective.
\end{lemma}

\noindent Next, the classical Riemann--Lebesgue lemma for Fourier transform states that: If $f\in L_1 (\R^n)$, then $(Ff)(y) \to 0$ as
$|y| \to \infty,$ and, hence, $Ff \in C_0 (\R^n).$ Therefore $\|Ff\|_{L_\infty (\R^n)} \leq \|f\|_{L_1 (\R^n)}$ and  $Ff$ is uniformly continuous. Based on this, we will demonstrate a similar version for the $\mathscr H$-transform as follows.

\begin{lemma}[Riemann--Lebesgue type for $\mathscr H$-transform]\label{riemann-lebegues} 
	If $f\in L_1(\mathbb R^n)$, then $\mathscr H f$ belongs to $C_0(\mathbb R^n)$; that is, $\mathscr H f$ is a bounded continuous function on $\R^n$, such that
	$
	\lim_{|y|\to\infty} (\mathscr H f)(y) = 0,
	$
	and satisfies the uniform estimate
	$
	\| \mathscr H f \|_{L_\infty(\mathbb R^n)} \le (2\pi)^{-n/2} \big(|a|+|b|\big)\,\|f\|_{L_1(\mathbb R^n)}.
	$
\end{lemma}
\begin{proof}
	For any function $f\in L_1(\mathbb R^n)$ the Fourier cosine and sine transforms given by \eqref{FsFc}
	are precisely the real and imaginary parts of the Fourier transform, that is
	$
	(Ff)(y)= (F_c f)(y)+ i (F_s f)(y).
	$ Therefore, we have
	$$(\mathscr H f)(y)= a(F_c f)(y) + b(F_s f)(y) = a.\textup{Re}\{(Ff)(y)\} + b.\textup{Im}\{(Ff)(y)\}.$$
	Since $|\cos(xy)| \le 1$ and $|\sin(x y)| \le 1$, we obtain
	$$
	|(\mathscr H f)(y)| \le (2\pi)^{-n/2}\,(|a|+|b|)\int_{\mathbb{R}^n} |f(x)|\,dx.
	$$
	Taking the supremum over $y$ yields, we claimed $L_\infty$ bound as	
	$$\| \mathscr H f \|_{L_\infty(\mathbb R^n)} \le (2\pi)^{-n/2} \big(|a|+|b|\big)\,\|f\|_{L_1(\mathbb R^n)}.
	$$
	Continuity of the $\mathscr H f$ directly follows from continuity of $Ff$ (by dominated convergence theorem) and linearity. Finally, by the Riemann–Lebesgue lemma for Fourier transform \cite{Rudin1991functional}, we infer $(Ff)(y)\to 0$ as $|y| \to \infty$, which implies that both $\textup{Re}\{(Ff)(y)\}$ and $\textup{Im}\{(Ff)(y)\}$ vanish at infinity; hence any linear combination of them vanishes too. Therefore $(\mathscr H f)(y)\to 0$ as $|y|$ tends to $\infty$. This completes the proof.
\end{proof}

\noindent Finally, we need the following lemma, which provides an estimate for the kernel function $K_{(a,b)}[f](x,v)$.
\begin{lemma}\label{lem2.1}
	Let $f \in L_q(\mathbb{R}^n)$ with $q \geq 1$. Then, for every $v \in \mathbb{R}^n$, the following estimate holds:
	\begin{equation}\label{eq2.3}
	\int_{\mathbb{R}^n} \left|K_{(a,b)}[f](x,v)\right|^q \, dx 
	\;\leq\; 4^{\,q-1}\Bigl(|3a^2-b^2|^q + 3(a^2+b^2)^q\Bigr)\,\|f\|_{L_q(\mathbb{R}^n)}^q.
	\end{equation}
\end{lemma}
\begin{proof}
	Recall the elementary inequality, valid for all complex numbers $\alpha_1,\dots,\alpha_k$ and $q \geq 1$,  $\big|\sum_{i=1}^k \alpha_i\big|^q 
	\;\leq\; k^{\,q-1}\sum_{i=1}^k |\alpha_i|^q,$
	which is a direct consequence of Hölder’s inequality (and reduces to the usual triangle inequality when $q=1$).
	From \eqref{eq2.2}, we observe that the kernel $K_{(a,b)}[f](x,v)$ consists of four terms. Applying the above inequality with $k=4$, we infer that
	$$\begin{aligned}
	&\int_{\mathbb{R}^n}\left|K_{(a,b)}[f](x,v)\right|^q dx
	\\&\leq \int_{\mathbb{R}^n}\bigg\{|(3a^2-b^2)f(x-v)| + (a^2+b^2)\big[\ |f(x+v)| + |f(-x+v)| + |f(-x-v)|\ \big]\bigg\}^q dx\\
	&\leq 4^{q-1} \bigg\{|3a^2-b^2|^q \int_{\mathbb{R}^n} |f(x-v)|^q dx + (a^2+b^2)^q \int_{\mathbb{R}^n} |f(x+v)|^q dx\\
	& +(a^2+b^2)^q \int_{\mathbb{R}^n} |f(-x+v)|^qdx + (a^2+b^2)^q \int_{\mathbb{R}^n} |f(-x-v)|^q dx\bigg\}.
	\end{aligned}$$
	
	\noindent 	Since $f \in L_q(\mathbb{R}^n)$, each of the above integrals equals $\|f\|_{L_q(\mathbb{R}^n)}^q$ by a suitable change of variables of the form $t = x \pm v$ and $t=-x \pm v$, with Jacobian determinant equal to $\pm 1$. Hence,
	$$
	\int_{\mathbb{R}^n} \left|K_{(a,b)}[f](x,v)\right|^q dx
	\;\leq\; 4^{\,q-1}\Bigl(|3a^2-b^2|^q + 3(a^2+b^2)^q\Bigr)\|f\|_{L_q(\mathbb{R}^n)}^q,
	$$
	which completes the proof.
\end{proof}
\section{Algebraic structures}\label{sec2}

\subsection{Existence, factorization and Parseval equalities}

\begin{theorem}\label{thm2.1}
	\noindent \textbf{(A)} 
	Let $f,g\in L_1(\mathbb{R}^n)$. Then the convolution \eqref{eq2.1} is well-defined for almost every $x\in\mathbb{R}^n$, belongs to $L_1(\mathbb{R}^n)$, and satisfies the norm inequality 
	\begin{equation}\label{eq2.4}
	\big\|f \underset{\mathscr{H}}{\ast} g\big\|_{L_1(\mathbb{R}^n)}
	\;\leq\; \frac{|3a^2-b^2| + 3(a^2+b^2)}{4|a|(2\pi)^{\tfrac{n}{2}}}
	\|f\|_{L_1(\mathbb{R}^n)} \, \|g\|_{L_1(\mathbb{R}^n)}.
	\end{equation}	 Moreover, the factorization identity holds
	
\begin{equation}\label{eq2.5}
	\mathscr{H}\big(f\underset{\mathscr{H}}{\ast} g\big)(y) = \big(\mathscr{H}f\big)(y) \big(\mathscr{H}g\big)(y) \ \ \text{for all}\ \ y\in\mathbb{R}^n. 	
\end{equation}		
	\noindent \textbf{(B)} 
	If $f,g\in L_2(\mathbb{R}^n)$, then $f\underset{\mathscr{H}}{\ast} g$ is defined for every $x$ and the factorization identity holds a.e. in $\mathbb{R}^n$. In addition, the following representation formula is valid:
	\begin{equation}\label{eq2.6}
	\big(f \underset{\mathscr{H}}{\ast} g\big)(x) 
	= \frac{1}{(2\pi)^{\tfrac{n}{2}}}\int_{\mathbb{R}^n} 
	\left(\frac{1}{a}\cos(xy) + \frac{1}{b}\sin(xy)\right)
	(\mathscr{H}f)(y)\,(\mathscr{H}g)(y)\, dy \ \ \text{a.e. in}\ \ \R^n.
	\end{equation}
	This formula can be regarded as a Parseval-type identity for the convolution structure \eqref{eq2.1}.	
\end{theorem}
\begin{proof}
	To prove the existence of the $\underset{\mathscr{H}}{\ast}$-convolution 
	 in $L_1(\mathbb{R}^n)$, it suffices to show that 
	$
	\int_{\mathbb{R}^{n}} \big|\big(f \underset{\mathscr{H}}{\ast} g\big)(x)\big|\, dx $ is finite.
	Indeed, from \eqref{eq2.1}, by a change of variables and an application of Fubini's theorem, we obtain
	\begin{align*}
	\int_{\mathbb{R}^{n}} \big|\big(f \underset{\mathscr{H}}{\ast} g\big)(x)\big|\, dx
	&\leq \frac{1}{4|a|(2\pi)^{\tfrac{n}{2}}} 
	\int_{\mathbb{R}^{2n}} |K_{(a,b)}[f](x,v)|\, |g(v)|\, dx\, dv  \\
	&= \frac{1}{4|a|(2\pi)^{\tfrac{n}{2}}} 
	\int_{\mathbb{R}^{n}} |g(v)| 
	\left(\int_{\mathbb{R}^n} |K_{(a,b)}[f](x,v)|\, dx\right) dv.
	\end{align*}
	By the estimation of the kernel function $K_{(a,b)}[f](x,v)$ in \eqref{eq2.3} with $q=1$, applied to $f\in L_1(\mathbb{R}^n)$, it follows that
	$
	\int_{\mathbb{R}^{n}} \big|\big(f \underset{\mathscr{H}}{\ast} g\big)(x)\big|\, dx 
	\leq \frac{|3a^2-b^2| + 3(a^2+b^2)}{4|a|(2\pi)^{\frac{n}{2}}}\,
	\|f\|_{L_1(\mathbb{R}^n)}\, \|g\|_{L_1(\mathbb{R}^n)}
	< +\infty.$
	This yields $f \underset{\mathscr{H}}{\ast} g \in L_1(\mathbb{R}^n)$, and thus the estimate \eqref{eq2.4} follows.
	To verify the factorization property, starting from the definition \eqref{eq1.3}, we compute
	\begin{align*}
	&(\mathscr{H}f)(y)(\mathscr{H}g)(y)
	\\&= \frac{1}{(2\pi)^n}\int_{\mathbb{R}^{2n}} (a\cos yu+b\sin yu)(a\cos yv+b\sin yv) f(u)g(v)\, du\, dv \\
	&= \frac{1}{4a(2\pi)^n} \int_{\mathbb{R}^{2n}} \Big\{(3a^2-b^2)\big[a\cos y(u+v) + b\sin y(u+v)\big] 
	+ (a^2+b^2)\big[a\cos y(u-v)\\&+ b\sin y(u-v) + a\cos y(-u+v)+b\sin y(-u+v)
	- a\cos y(-u-v) - b\sin y(-u-v)\big]\Big\} f(u)g(v)\, du\, dv.
	\end{align*}
	After rearranging and applying the substitutions 
	$t = u\pm v$ and $t=-u\pm v$, we deduce  
	\begin{align*}
(\mathscr{H}f)(y)(\mathscr{H}g)(y)
	&= \frac{1}{4a(2\pi)^n} \bigg\{\int_{\mathbb{R}^{2n}}(3a^2-b^2)(a\cos yt+b\sin yt) f(t-v)g(v)\, dt\, dv \\&+(a^2+b^2)\Big[\int_{\mathbb{R}^{2n}}(a\cos yt+b\sin yt) f(t+v) g(v)\, dt\, dv 
	\\&+\int_{\mathbb{R}^{2n}}(a\cos yt+b\sin yt) f(-t+v) g(v)\, dt\, dv -\int_{\mathbb{R}^{2n}}(a\cos yt+b\sin yt) f(-t-v) g(v)\, dt\, dv\Big]\bigg\} \\
	&= \frac{1}{(2\pi)^{\frac{n}{2}}} \frac{1}{4a(2\pi)^{\frac{n}{2}}}
	\int_{\mathbb{R}^n}(a\cos yt+b\sin yt) 
	\bigg(\int_{\mathbb{R}^n}\Big[(3a^2-b^2)f(t-v) + (a^2+b^2)\big(f(t+v)\\&+f(-t+v)-f(-t-v)\big)\Big] g(v)\, dv\bigg) dt.
	\end{align*}
	Combining \eqref{eq2.1} and \eqref{eq2.2}, this expression simplifies to
	$$
	(\mathscr{H}f)(y)(\mathscr{H}g)(y)
	= \frac{1}{(2\pi)^{\frac{n}{2}}}\int_{\mathbb{R}^n} (a\cos yt+b\sin yt)\, \big(f \underset{\mathscr{H}}{\ast} g\big)(t)\, dt
	= \mathscr{H}\big(f \underset{\mathscr{H}}{\ast} g\big)(y).
	$$
	Since, $f,g \in L_1 (\R^n),$ then $\mathscr{H}f$ and $\mathscr{H}g$ are continuous, therefore, the identity holds pointwise, that proves \eqref{eq2.5}.
	Finally, when $f,g\in L_2(\mathbb{R}^n)$, an analogous computation yields the representation formula \eqref{eq2.6}. The proof is complete.
\end{proof}

\begin{remark}[Derivation and naturality of the kernel function]
	\textup{The specific form of the convolution kernel $K_{(a,b)}$ in the formula \eqref{eq2.2} is not arbitrary but arises naturally from the fundamental requirement that the $\mathscr{H}$-transform should be an algebra homomorphism, i.e., $\mathscr{H}\big(f\underset{\mathscr{H}}{\ast} g\big)(y) = \big(\mathscr{H}f\big)(y) \big(\mathscr{H}g\big)(y)$, which is also the factorization identity \eqref{eq2.5} has just been proven above. To see why kernel in \eqref{eq2.2} is the unique (up to normalization) kernel satisfying \eqref{eq2.5}, we expand the right-hand side using the definition \eqref{eq1.3} as follows:
		\begin{align*}
		(\mathscr{H}f)(y)(\mathscr{H}g)(y) 
		= \frac{1}{(2\pi)^{n}} \iint_{\mathbb{R}^{2n}} \bigl(a\cos(xy)+b\sin(xy)\bigr)
		\times \bigl(a\cos(vy)+b\sin(vy)\bigr) f(x)g(v)\,dx\,dv.
		\end{align*}
		Applying elementary trigonometric identities,
		\begin{align*}
		\cos(xy)\cos(vy) &= \tfrac12\bigl[\cos((x-v)y) + \cos((x+v)y)\bigr],\\
		\sin(xy)\sin(vy) &= \tfrac12\bigl[-\cos((x-v)y) + \cos((x+v)y)\bigr],\\
		\cos(xy)\sin(vy) &= \tfrac12\bigl[\sin((x+v)y) - \sin((x-v)y)\bigr],\\
		\sin(xy)\cos(vy) &= \tfrac12\bigl[\sin((x+v)y) + \sin((x-v)y)\bigr],
		\end{align*}
		we obtain after simplification
		\begin{equation}\label{hfhg}
		\begin{aligned}
		&(\mathscr{H}f)(y)(\mathscr{H}g)(y) 
		\\&= \frac{1}{(2\pi)^{n}} \iint_{\mathbb{R}^{2n}} \Bigl[
		\frac{a^{2}-b^{2}}{2}\cos((x-v)y) 
		+ \frac{a^{2}+b^{2}}{2}\cos((x+v)y) + ab\sin((x+v)y) \Bigr] f(x)g(v)\,dx\,dv.
		\end{aligned}
		\end{equation}
		The left-hand side of \eqref{eq2.5} is, by definition of the convolution, we have
		$$
		\mathscr{H}(f\underset{\mathscr{H}}{\ast}g)(y) = \frac{1}{(2\pi)^{n/2}} \int_{\mathbb{R}^{n}} \bigl[a\cos(xy)+b\sin(xy)\bigr](f\underset{\mathscr{H}}{\ast}g)(x)\,dx.
		$$
		Inserting the ansatz $$
		(f\underset{\mathscr{H}}{\ast}g)(x)=\frac{1}{4a(2\pi)^{\frac{n}{2}}}\int_{\mathbb{R}^{n}} K_{(a,b)}[f](x,v)\,g(v)\,dv,
		$$
		with
		$
		K_{(a,b)}[f](x,v)= C_{1}f(x-v)+C_{2}f(x+v)+C_{3}f(-x+v)+C_{4}f(-x-v),
		$
		and carrying out the integration, we infer an expression of the form
		$$
		\mathscr{H}(f\underset{\mathscr{H}}{\ast}g)(y)=\frac{1}{(2\pi)^{n}} \iint_{\mathbb{R}^{2n}} \bigl(a\cos(xy)+b\sin(xy)\bigr)
		\bigl[C_{1}f(x-v)+\cdots +C_{4}f(-x-v)\bigr]g(v)\,dx\,dv.
		$$
		After the changes of variables $t=x\pm v$, $t=-x\pm v$ and comparing the resulting coefficients of the trigonometric terms with those in \eqref{hfhg}, we are led to a linear system for $C_{1},\dots ,C_{4}$ after simplification:
		$$\left\{\begin{aligned}
		a C_1+a C_3 & =\frac{a^2-b^2}{2} \\
		a C_2+a C_4 & =\frac{a^2+b^2}{2} \\
		b C_2-b C_4 & =a b
		\end{aligned}\right.$$
		Solving this system gives the unique solution: $C_{1}= \frac{3a^{2}-b^{2}}{4a}$, $C_{2}=C_{3}= \frac{a^{2}+b^{2}}{4a}$ and, $C_{4}= -\frac{a^{2}+b^{2}}{4a}$. Multiplying these coefficients by $4a$ (to absorb the factor $\frac{1}{4a}$ in front of the integral in \eqref{eq2.1}) yields precisely the kernel displayed in \eqref{eq2.2}.  Consequently, the kernel $K_{(a,b)}$ is not chosen ad hoc; it is the unique bilinear expression in the shifts of $f$ that makes the diagram}
	\begin{center}
		\begin{tikzcd}[column sep=1.0cm]
		L_{1}(\mathbb{R}^{n})\times L_{1}(\mathbb{R}^{n}) \arrow[r,"{\underset{\mathscr{H}}{\ast}}"] \arrow[d,"{\mathscr{H}\times\mathscr{H}}"'] &
		L_{1}(\mathbb{R}^{n}) \arrow[d,"\mathscr{H}"] \\
		C_{0}(\mathbb{R}^{n})\times C_{0}(\mathbb{R}^{n}) \arrow[r,"\textup{pointwise}",
		"\textup{product}"']&
		C_{0}(\mathbb{R}^{n})
		\end{tikzcd}
	\end{center}
	\textup{commute. This derivation underscores the naturality of the $\underset{\mathscr{H}}{\ast}$-convolution product: it is forced by the algebraic requirement that $\mathscr{H}$-transform \eqref{eq1.3} be a homomorphism.}
\end{remark}

\subsection{Main theorems}\label{sub3.2}
\begin{theorem}[Banach algebra]\label{theorem22}
If $\underset{\mathscr{H}}{\ast}$ is regarded as the multiplication on $L_1(\mathbb{R}^n)$, then the structure  $(L_1(\mathbb{R}^n),+,\underset{\mathscr{H}}{\ast})$ equipped with the norm defined by $\|f\|_{\a,L_1(\mathbb{R}^n)}:=\int_{\mathbb{R}^n} \alpha|f(x)| dx$ for all $f\in L_1(\mathbb{R}^n)$, where $ \alpha=\sqrt{\frac{|3a^2-b^2|+3(a^2+b^2)}{4|a|(2\pi)^{\frac{n}{2}}}}>0$, forms a commutative Banach algebra over $\mathbb{C}$ without an identity element.
\end{theorem}
\begin{proof}
Denote $\mathcal{A}=\big(L_1(\mathbb{R}^n),+,\underset{\mathscr{H}}{\ast},\|.\|_{\a,L_1}\big)$. We verify that $\mathcal{A}$ satisfies all axioms of a commutative Banach algebra, and then prove the nonexistence of an identity element of multiplicative. This proof is divided into three steps. \\ \textbf{Step 1: Algebraic axioms.}  By Theorem~\ref{thm2.1}, the space $L_1(\mathbb{R}^n)$ is closed under the convolution \eqref{eq2.1}. Since $\underset{\mathscr{H}}{\ast}$ is defined via the $\mathscr{H}$-transform, and $L_1(\mathbb{R}^n)$ is already a complex vector space under addition and scalar multiplication, it follows that $\mathcal{A}$ is a complex vector space.  We next verify bilinearity and distributivity.  Indeed, let $f,g,h \in L_1(\mathbb{R}^n)$ and $\lambda \in \C$. Utilizing the linearity of $\mathscr{H}$-transform \eqref{eq1.3} and factorization property \eqref{eq2.5} yields
\[
\begin{aligned}
\mathscr{H}\big[(f + g) \underset{\mathscr{H}}{\ast} h\big](y) &= \mathscr{H}(f + g)(y). (\mathscr{H} h)(y) \\
&= \big[(\mathscr{H} f)(y) + (\mathscr{H} g)(y)\big] .(\mathscr{H} h)(y) \\
&= (\mathscr{H} f)(y)(\mathscr{H} h)(y) + (\mathscr{H} g)(y)(\mathscr{H} h)(y) \\
&= \mathscr{H}(f \underset{\mathscr{H}}{\ast} h)(y) + \mathscr{H}(g \underset{\mathscr{H}}{\ast} h)(y) 
= \mathscr{H}\big[(f \underset{\mathscr{H}}{\ast} h) + (g \underset{\mathscr{H}}{\ast} h)\big](y).
\end{aligned}
\]
By Lemma \ref{donanh}, the injectivity of $\mathscr{H}$-transform \eqref{eq1.3} is valid. Therefore, it follows that $(f+g) \underset{\mathscr{H}}{\ast} h = (f \underset{\mathscr{H}}{\ast} h) + (g \underset{\mathscr{H}}{\ast} h),$ which establishes right distributivity. Left distributivity follows by an analogous argument, meaning  $f \underset{\mathscr{H}}{\ast} (g + h) = (f \underset{\mathscr{H}}{\ast} g) + (f \underset{\mathscr{H}}{\ast} h).$ Moreover, since 
$$\begin{aligned}
	\mathscr{H}[\lambda(f \underset{\mathscr{H}}{\ast} g)](y) & =\lambda \mathscr{H}(f \underset{\mathscr{H}}{\ast} g)(y)=\lambda(\mathscr{H} f)(y)(\mathscr{H} g)(y) \\
	& =(\mathscr{H}(\lambda f))(y)(\mathscr{H} g)(y)=\mathscr{H}[(\lambda f)\underset{\mathscr{H}}{\ast} g](y) .
\end{aligned}$$
By the injectivity of $\mathscr{H}$-transform, then $\lambda(f \underset{\mathscr{H}}{\ast} g)=(\lambda f)\underset{\mathscr{H}}{\ast} g$. A symmetric argument shows $\lambda(f \underset{\mathscr{H}}{\ast} g)= f\underset{\mathscr{H}}{\ast} (\lambda g)$. Therefore, homogeneity for scalar multiplication is also satisfied, and the bilinearity of $\underset{\mathscr{H}}{\ast}$ is established.
Again, applying the factorization property \eqref{eq2.5} repeatedly, we achieve
\begin{align*}
\mathscr{H}\big[f \underset{\mathscr{H}}{\ast}\big(g \underset{\mathscr{H}}{\ast} h\big)\big](y) &= (\mathscr{H}f)(y) \mathscr{H}\big(g \underset{\mathscr{H}}{\ast} h\big)(y) \\&= (\mathscr{H}f)(y)(\mathscr{H}g)(y) (\mathscr{H}h)(y)\\
&=\mathscr{H}(f \underset{\mathscr{H}}{\ast} g)(y). (\mathscr{H} h)(y) = \mathscr{H}\big[\big(f \underset{\mathscr{H}}{\ast} g\big)\underset{\mathscr{H}}{\ast} h\big](y).
\end{align*}Hence, $f\underset{\mathscr{H}}{\ast}(g \underset{\mathscr{H}}{\ast} h)  =(f \underset{\mathscr{H}}{\ast} g) \underset{\mathscr{H}}{\ast} h,$ which implies that the associativity is valid. Finally, by \eqref{eq2.5} and the commutativity of multiplication in $\C$, we have 
$\mathscr{H}(f \underset{\mathscr{H}}{\ast}g)(y)=(\mathscr{H}f)(y)(\mathscr{H}g)(y)=(\mathscr{H}g)(y)(\mathscr{H}f)(y)=\mathscr{H}(g\underset{\mathscr{H}}{\ast}f)(y).$ This means that $f \underset{\mathscr{H}}{\ast}g= g \underset{\mathscr{H}}{\ast}f$.\\
\textbf{Step 2: Analytic axioms.} We now show that the norm is submultiplicative and that the space is complete. 
Indeed, the norm $\|.\|_{\a,L_1}$ is a scaled version of the standard $L_1$-norm, so it is indeed a norm. Its completeness follows from the completeness of $(L_1 (\R^n), \|.\|_{L_1})$ and the fact that the norms $\|.\|_{\a,L_1}$ and $\|.\|_{L_1}$ are equivalent.
The crucial submultiplicative property is provided directly by the inequality \eqref{eq2.4} as follows
$$
\big\|f \underset{\mathscr{H}}{\ast} g\big\|_{\a,L_1(\mathbb{R}^n)} \leq \bigg(\, \int_{\mathbb{R}^n} \alpha |f(x)| dx\bigg)\bigg(\int_{\mathbb{R}^n} \alpha|g(x)| dx\bigg)= \|f\|_{\a,L_1(\mathbb{R}^n)} \|g\|_{\a,L_1(\mathbb{R}^n)}, 
$$ where $\alpha=\sqrt{\frac{|3a^2-b^2|+3(a^2+b^2)}{4|a|(2\pi)^{\frac{n}{2}}}}.$
This confirms that $\mathcal{A}$ is a complete normed algebra (normed ring \cite{1964}). All of these lead to the assertion that $\mathcal{A}$ is a Banach algebra.\\
\textbf{Step 3: Nonexistence of an identity element.}
We conclude by proving that 
$\mathcal{A}$ lacks multiplicative identity element.
Assume, for contradiction, that there exists $e\in L_1(\mathbb{R}^n)$ such that
$
e\underset{\mathscr{H}}{\ast}f=f,$ for any functions $f\in L_1(\mathbb{R}^n).
$
Applying $\mathscr{H}$-transform \eqref{eq1.3} to both sides and using \eqref{eq2.5},   we get
$
(\mathscr{H} e)(y)(\mathscr{H} f)(y)=(\mathscr{H} f)(y), \forall  y \in \mathbb{R}^n. $
This implies
$(\mathscr{H} f)(y)\big[(\mathscr{H} e)(y)-1\big]=0.$
Since the image of $\mathscr{H}$-transform \eqref{eq1.3} on $L_1\left(\mathbb{R}^n\right)$ is sufficiently rich (i.e., it contains functions that are nonzero almost everywhere), we must have
$
(\mathscr{H} e)(y)=1$ for almost all  $y \in \mathbb{R}^n.
$
However, this leads to a contradiction. By Lemma \ref{riemann-lebegues}, for any $\mathscr{H}e$ with $e\in L_1(\mathbb{R}^n)$ must vanish at infinity, that is
$
\lim_{|y|\to\infty} (\mathscr{H}e)(y)=0.
$
This contradicts the fact that $(\mathscr{H}e)(y)=1$ almost everywhere. Thus, no such $e$ exists, and $\mathcal{A}$ has no multiplicative identity element. This completes the proof. 
\end{proof}
Theorem \ref{theorem22} provides an affirmative answer to the question raised in Section \ref{sec1} and also constitutes the first main contribution of this paper.
The next following theorem shows the absence of zero divisors for the convolution \eqref{eq2.1}. It can be regarded as an analogue of Titchmarsh's classical convolution theorem, adapted to the present framework.

\begin{theorem}[Titchmarsh type theorem]\label{Titchmarsh}	Let $f,g\in L_1(\mathbb{R}^n)$ be functions with compact support. 
	If
	$
	(f\underset{\mathscr H}{\ast} g)(x)=0$ for almost every $x\in\mathbb{R}^n$,
	then either $f=0$ almost everywhere or $g=0$ almost everywhere on $\mathbb{R}^n$.
	
\end{theorem}
\begin{proof}
	By the Paley--Wiener theorem (see {\bf 7.22}, p.198 in \cite{Rudin1991functional}), if $f$ has a compact support, then the Fourier transform $Ff(z)$ extends to an entire function on $\mathbb{C}^n$; the same is true for the inverse transform (i.e. $F^{-1}f(z)$). Since $Ff$ and $F^{-1}f$ extend to entire functions on $\mathbb{C}^n$, then the linear combination $\mathscr{H}f(z) = \frac{a-ib}{2}Ff(z) + \frac{a+ib}{2}F^{-1}f(z),
	$ also extends to an entire function; consequently, its restriction to $\mathbb{R}^n$ is real‑analytic. The same holds for $\mathscr{H}g(z)$. For a real‑analytic function that is not identically zero, its zero set has empty interior. This is the only property needed in the subsequent Baire‑category argument.
	Hence, the restrictions $\mathscr{H}f|_{\mathbb{R}^n}$ and $\mathscr{H}g|_{\mathbb{R}^n}$ are real-analytic functions on $\mathbb{R}^n$, whose zero sets cannot contain open subsets unless the function is identically zero.
	From the factorization property \eqref{eq2.5} $\forall y\in\mathbb{R}^n$, we get
	$$
	\mathscr{H}\big(f\underset{\mathscr H}{\ast} g\big)(y) = (\mathscr{H}f)(y)\,(\mathscr{H}g)(y).
	$$
	Because of the assumption $f\underset{\mathscr H}{\ast} g=0$ almost every on $\mathbb{R}^n$, then the left-hand side vanishes for all $y\in\mathbb{R}^n$. Thus
	$
	(\mathscr{H}f)(y)\,(\mathscr{H}g)(y)=0$  for all $y\in\mathbb{R}^n.
	$
	Now we denote zero-sets
	$$Z_f:=\big\{y\in\mathbb{R}^n :(\mathscr{H}f)(y)=0\big\} \quad \text{and} \quad 
	Z_g:=\big\{y\in\mathbb{R}^n :(\mathscr{H}g)(y)=0\big\},$$ then $Z_f\cup Z_g=\mathbb{R}^n$. If one of $\mathscr{H}f$ or $\mathscr{H}g$ is identically zero on $\mathbb{R}^n$, we are done; suppose instead that neither is identically zero. Then each of $\mathscr{H}f$ and $\mathscr{H}g$ is a nontrivial real-analytic function on the path-connected $\mathbb{R}^n$. A standard property of nontrivial real-analytic functions is that their zero sets have empty interior (equivalently, they are nowhere dense in $\mathbb{R}^n$). Hence $Z_f$ and $Z_g$ are closed and nowhere dense subsets of $\mathbb{R}^n$.	However, $\mathbb{R}^n$ (with its usual metric) is a Baire space; the Baire category theorem (see {\bf 5.6}, p.97 in \cite{Rudin87}, also \cite{Stein2011}) asserts that a complete metric space cannot be expressed as the union of two nowhere dense sets. This contradicts $Z_f\cup Z_g=\mathbb{R}^n$. Therefore, our assumption that both $\mathscr{H}f$ and $\mathscr{H}g$ are nontrivial is false, which implies that at least one of them vanishes identically on $\mathbb{R}^n$. Without loss of generality, we assume $\mathscr{H}f\equiv 0$ on  $\mathbb{R}^n.$ Since \eqref{eq1.4} gives an explicit inverse operator $\mathscr H^{-1}$ on the class under consideration, applying the inverse operator $\mathscr H^{-1}$ on the both sides, we conclude that $f=\mathscr H^{-1}(0)=0$ almost everywhere on $\R^n$. The alternative case $\mathscr{H}g\equiv 0$ is analogous. This completes the proof.
\end{proof}
\begin{proposition}[Density property]\label{density-H-image}
	The image of compactly supported smooth functions under $\mathscr H$-transform \eqref{eq1.3} i.e.,
	$
	\mathscr H(C_c^\infty(\mathbb R^n)),
	$
	is dense in $C_0(\mathbb R^n)$ space with respect to the uniform norm. Consequently $\mathscr H(L_1(\mathbb R^n))$ is dense in $C_0(\mathbb R^n)$.
\end{proposition}
\begin{proof}
	We divide the proof into three steps.\\	
	\textbf{Step 1: Algebra generated by characters is dense on each compact set.}
	For each fixed $x\in\mathbb R^n$, we define the character $\chi_x:\mathbb R^n\to\mathbb C,$ by $\chi_x(y):=e^{i x\cdot y}.$ Consider the complex algebra
	$
	\mathcal B := \operatorname{span}\big\{\chi_x : x\in\mathbb R^n\big\}
	$
	of finite linear combinations of characters. Fix a compact set $K\subset\mathbb R^n$ and view $\mathcal B$ as an algebra of continuous functions on $K$. We verify the conditions of the Stone--Weierstrass theorem (see Chapter 5 in \cite{Rudin1991functional}) on compact $K \subset \R^n$ as follows:
	
\noindent	(i) $\mathcal B$ is an algebra containing the constant functions (take $x=0$, then $\chi_0\equiv 1$).
	
\noindent	(ii) $\mathcal B$ separates points of $K$. Indeed, if $y_1,y_2\in K$ with $y_1\neq y_2$. Put $v:=y_1-y_2$ and choose $\lambda\in\mathbb R$ so that $\lambda |v|^2\notin 2\pi\mathbb Z$ (certainly feasible to choose
$\lambda$, because
$\R$ is dense in
$\R/2\pi\mathbb{Z}$). Then 
$
	\chi_{\lambda v}\left(y_1\right)=e^{i \lambda v \cdot y_1}$ and $  \chi_{\lambda v}\left(y_2\right)=e^{i \lambda v \cdot y_2}.$
	Since $\lambda v \left(y_1-y_2\right)=\lambda|v|^2 \notin 2 \pi \mathbb{Z}$, we infer $\chi_{\lambda v}\left(y_1\right) \neq \chi_{\lambda v}\left(y_2\right),$
 hence some character separates $y_1$ and $y_2$.
	
\noindent	(iii) $\mathcal B$ is closed under complex conjugation, since $\overline{\chi_x}=\chi_{-x}\in\mathcal B$.
	
\noindent By Stone--Weierstrass's theorem (applied to the compact space $K$), we infer the algebra $\mathcal B$ is uniformly dense in $C(K)$.
	
\noindent \textbf{Step 2: Approximation by Fourier transform of smooth compactly supported functions.}
	Let $\varphi\in C(K)$ and $\varepsilon>0$. By \textbf{Step 1}, there exists a finite linear combination of characters (trigonometric polynomial)
$
	p(y)=\sum_{j=1}^N c_j e^{i x_j\cdot y},$ where $x_j\in\mathbb R^n$ and $c_j\in\mathbb C,
$
	such that \begin{equation}\label{111}
	\sup_{y\in K}|p(y)-\varphi(y)|<\varepsilon/2.
	\end{equation}
	\noindent	Each exponential $e^{i x_j\cdot y}$ is the Fourier transform (in the variable \(y\)) of the Dirac measure at \(x_j\). To obtain an approximation by the Fourier transform of a compactly supported smooth function, choose a standard mollifier $\rho\in C_c^\infty(\mathbb R^n)$ with $\int\rho=1$ and put $\rho_\delta(x)=\delta^{-n}\rho(x/\delta)$. Define
	$
	\varphi_{j,\delta}(y):=F(\rho_\delta * \delta_{x_j})(y)=e^{i x_j\cdot y}\,(F\rho_\delta)(y).
	$ Since $(F\rho_\delta)(y)\to 1$ uniformly on compact K as $\delta\to 0$, for sufficiently small $\delta$ we have
	$$\sup _{y \in K}|e^{i x_j \cdot y}-\varphi_{j, \delta}(y)|=\sup _{y \in K}\big|e^{i x_j \cdot y}\left[1-(F\rho_\delta)(y)\right]\big| \leq \sup _{y \in K}|1-(F\rho_\delta)(y)|<\frac{\varepsilon}{2 N \max _j\left|c_j\right|}.$$
Set $p_\delta(y):=\sum_{j=1}^N c_j\varphi_{j,\delta}(y),
	$
then	$p_\delta$ is the Fourier transform of a function in $C_c^\infty(\mathbb R^n)$  and
	\begin{equation}\label{222}
	\left|p(y)-p_\delta(y)\right|=\bigg|\sum_{j=1}^N c_j\big(e^{i x_j \cdot y}-\varphi_{j, \delta}(y)\big)\bigg| \leq \sum_{j=1}^N\left|c_j\right|\ |e^{i x_j \cdot y}-\varphi_{j, \delta}(y)|<\sum_{j=1}^N\left|c_j\right| \frac{\varepsilon}{2 N \max _j\left|c_j\right|} \leq \varepsilon / 2 .
	\end{equation}
\noindent	Coupling \eqref{111} and \eqref{222} yields
	$$
	\sup _{y \in K}\left|p_\delta(y)-\varphi(y)\right| \leq \sup _{y \in K}\left|p_\delta(y)-p(y)\right|+\sup _{y \in K}|p(y)-\varphi(y)|<\varepsilon / 2+\varepsilon / 2=\varepsilon.
	$$
	Hence, every continuous function on the compact $K$ can be uniformly approximated on $K$ by Fourier transforms of compactly supported smooth functions.\\
\textbf{Step 3: Global approximation and passage to $\mathscr H$-transform.}
Let $g\in C_0(\mathbb R^n)$ and fix $\varepsilon>0$. Since $g$ vanishes at infinity, there exists a compact $K$ with $\sup_{y\not\in K}|g(y)|<\varepsilon/3$. By \textbf{Step 2}, there exists a function $\Phi=F\psi$  with $\psi\in C_c^\infty$ such that
	$
	\sup_{y\in K}|g(y)-\Phi(y)|<\varepsilon/3.
$
Moreover, since both $g$ and $\Phi$ are bounded functions, we deduce
$$
	\sup_{y\in\mathbb R^n}|g(y)-\Phi(y)| \le \sup_{y\in K}|g-\Phi| + \sup_{y\notin K}|g| + \sup_{y\notin K}|\Phi| .
$$
Since $\Phi=F\psi$ with $\psi\in C_c^\infty$, then $\Phi \in C_0 (\R^n)$ due to classical Riemann--Lebesgue lemma \cite{Rudin1991functional}. Hence, by choosing the compact $K$ large enough and the mollifier parameter small enough, we may also ensure $\sup_{y\notin K}|\Phi(y)|<\varepsilon/3.$ Combining the above arguments, we get $\sup_{y\in\mathbb R^n}|g(y)-\Phi(y)|<\varepsilon$. This shows that the set $\big\{F\psi:\psi\in C_c^\infty(\mathbb R^n)\big\}$ is dense in $C_0(\mathbb R^n)$ in the uniform norm.	
	Recall that $\mathscr H$-transform \eqref{eq1.3} is  a fixed real linear combination of the Fourier cosine and sine transforms, and the latter provide the real and imaginary parts of the Fourier transform, i.e, from \eqref{eq1.1} we have 
	$
	(a-i b) F f(y)=(2\pi)^{-n/2} \int_{\R^n}(a-i b) e^{i x \cdot y} f(x) d x,
	$
	then
	$\textup{Re}\big\{(a-i b) e^{i x \cdot y}\big\}=a \cos (xy)+b \sin (xy),
	$ leading to
	$
	\mathscr{H} f=\textup{Re}\{(a-i b) F f\}.
	$	
This means that $\mathscr H(C_c^\infty)$ is exactly the real linear image of set $\{F\psi:\psi\in C_c^\infty\}$ under a continuous linear mapping $T: C_0(\mathbb R^n;\mathbb C)\to C_0(\mathbb R^n;\mathbb R)$. By definition in \eqref{eq1.3} with the condition $a \neq 0$ and $b \neq 0$, then $a-i b \neq 0$, which yields  $T$ is a continuous surjective linear map. 
Therefore $\mathscr H(C_c^\infty(\mathbb R^n))$ is dense in $C_0(\mathbb R^n)$, and the same holds for $\mathscr H(L_1 (\mathbb R^n))$ since $C_c^\infty (\mathbb R^n)\subset L_1 (\mathbb R^n)$ and $\mathscr H:L_1 (\mathbb R^n)\to C_0 (\mathbb R^n)$ is continuous by Lemma \ref{riemann-lebegues}.
\end{proof}

Proposition \ref{density-H-image} is a useful tool for proving the following result for $\mathscr H$-algebras, also known as the Wiener–Lévy invertibility criterion.

\begin{theorem}[Wiener--Lévy invertibility criterion for $\mathscr H$-algebras]\label{WienerHalgebra}
Let $\mathscr H: L_1 (\R^n) \to C_0 (\R^n)$	be the transform is defined by \eqref{eq1.3}. Equip $L_1 (\R^n)$ with the product $f \cdot g:=f\underset{\mathscr H}{\ast} g$ is the convolution \eqref{eq2.1}. Let $g \in L_1 (\R^n)$ and suppose that
$$
1+(\mathscr H g)(y)\neq 0 \ \text{for all }y\in\mathbb R^n.$$
Then there exists $\ell\in L_1(\mathbb R^n)$ such that for every $y \in \R^n$	
	$$
	(\mathscr H \ell)(y)=-\frac{(\mathscr H g)(y)}{1+(\mathscr H g)(y)}.
	$$
Equivalently, the unitized element 
$(g,1)$ is invertible in the unitization of the 
$\mathscr H$-algebras, and the transform of its inverse yields the displayed rational function.
\end{theorem}

\begin{proof}
Let $\mathcal{A}=\big(L_1(\mathbb{R}^n),\underset{\mathscr{H}}{\ast},\|.\|_{\a,L_1(\mathbb{R}^n)}\big)$ denote $L_1 (\R^n)$ equipped with the $\underset{\mathscr{H}}{\ast}$-convolution; (note $\|.\|_{\a,L_1}$ and $\|.\|_{L_1}$ are equivalent). By Theorem \ref{theorem22}, then $\mathcal{A}$ is a commutative Banach algebra (without unit).
Defined $$\Gamma:\mathcal{A}\to C_0 (\R^n)\ \ \text{by}\ \ \Gamma(f)=\mathscr{H} f.$$ The multiplicativity of the $\mathscr H$-transform shows that $\Gamma$ is a continuous algebra homomorphism from $\mathcal{A}$ into the commutative Banach algebra $C_0 (\R^n)$, which implies 
$
\Gamma(f\cdot g)=\Gamma(f)\Gamma(g).
$
Moreover, by Lemma \ref{riemann-lebegues} then the map $\Gamma$ is continuous and	$\|\Gamma \|_{L_\infty(\mathbb R^n)}=\| \mathscr H f \|_{L_\infty(\mathbb R^n)} \le (2\pi)^{-n/2} \big(|a|+|b|\big)\,\|f\|_{L_1(\mathbb R^n)}.
$
Let $\mathfrak{M}(\mathcal{A})$ denote the maximal ideal space (the set of nonzero multiplicative linear functionals on $\mathcal{A}$).
First, we need to identify the maximal ideal space. For each  fixed $y\in\mathbb R^n$, the evaluation $\varphi_y(f):=(\mathscr{H} f)(y)$ is a nonzero multiplicative linear functional on $\mathcal{A}$. We will show that every character of $\mathcal{A}$ has this form. Indeed, 
define $$\psi:\Gamma(\mathcal{A})\to\mathbb C\ \ \text{by}\ \ \psi(\Gamma f)=\varphi(f).$$ 
Since $\Gamma$ is an algebra homomorphism and $\varphi$ is multiplicative, $\psi$ is well-defined and multiplicative on $\Gamma(\mathcal{A})$. By Proposition \ref{density-H-image}, then $\Gamma(\mathcal{A})$ is uniformly dense in $C_0(\mathbb R^n)$. Moreover,  the subalgebra $\Gamma(\mathcal{A}) \subset C_0(\mathbb R^n)$ separates points (because for $y_1\neq y_2$, one can find $f$ with $\mathscr H f(y_1)\neq\mathscr H f(y_2)$ using Urysohn's lemma; similar argument {\bf 11.13}, p.283 in \cite{Rudin1991functional} and density property) and contains the constants (take a suitable Gaussian $f$ with $\mathscr H f\equiv 1$). Hence the uniform closure of $\Gamma(\mathcal{A})$ is $C_0(\mathbb R^n)$ itself. 
On the other hand, we have inequality
\[
|\psi(\Gamma f)| = |\varphi(f)| \le \|\varphi\|\|f\|_{\alpha} \le M.\|\Gamma f\|_\infty,\ \text{where}\ M \ \text{is a positive constant},
\] then  the functional $\psi$ is continuous on $\Gamma(\mathcal{A})$ and extends uniquely to a continuous multiplicative linear functional on $C_0(\mathbb R^n)$. Every such functional on $C_0(\mathbb R^n)$ is evaluation at some point $y\in\mathbb R^n$. Thus there exists $y\in\mathbb R^n$ such that $\psi(h)=h(y)$ for all $h\in C_0(\mathbb R^n)$. Consequently,
\[
\varphi(f)=\psi(\Gamma f)=(\mathscr H f)(y)=\varphi_y(f).
\]
Therefore $$\mathfrak{M}(\mathcal{A})=\{\varphi_y:y\in\mathbb R^n\}.$$

\noindent Under this identification, Gelfand's transform (see {\bf 2.2.4}, p.54 in \cite{Kaniuth09}) is define $\hat{f}: \mathcal{A} \rightarrow C(\mathfrak{M}(\mathcal{A}))$ satisfies
$$
\hat{f}\left(\varphi_y\right)=\varphi_y(f)=(\mathscr{H} f)(y),
$$
so $\hat{f}$ coincides with $\mathscr{H} f$ regarded as a continuous function on $\mathbb{R}^n$. We now handle the unitarization using the Gelfand invertibility criterion, which is also key to this prove. We set
 $\mathcal A_1=\mathcal A\oplus\mathbb C$ with product
	\[
	(f,\beta)\cdot(g,\gamma)=(f\cdot g+\beta g+\gamma f,\ \beta\gamma),
	\]
	and norm $\|(f,\beta)\|'=\textup{Const}.\|f\|_{L_1 (\R^n)}+|\beta|$. Then $\mathcal A_1$ is also a unital commutative Banach algebra with unit element $(0,1)$. Its maximal ideal space may be identified with that of $\mathcal A$, and the Gelfand transform of $(f,\beta)$ is
	\[
	\widehat{(f,\beta)}(\varphi_y)=(\mathscr H f)(y)+\beta,\quad y\in\mathbb R^n.
	\] 
By the Gelfand invertibility criterion (see \cite{Kaniuth09}) for unital commutative Banach algebras in $\mathcal{A}_1$: an element of $\mathcal{A}_1$ is invertible if and only if its Gelfand transform vanishes nowhere on the maximal ideal. Consider $(g, 1) \in \mathcal{A}_1$. Its Gelfand transform is
$
\widehat{(g,1)}(\varphi_y)=(\mathscr H g)(y)+1, y\in\mathbb R^n.
$ 
By assumption $1+(\mathscr H g)(y)\neq0$ for all $y \in \R^n$, precisely guarantees this non-vanishing. Hence $\widehat{(g,1)}(y)\neq0$ on $\mathfrak M(\mathcal A)$, which implies that $(g,1)$ is invertible in $\mathcal A_1$. 
Thus there exists  $(\ell,\lambda)\in\mathcal A_1$ such that $$ (g,1)\cdot(\ell,\lambda)=(0,1). $$
Expanding the left-hand products yields
$
(g, 1) \cdot(\ell, \lambda)=(g \cdot \ell+\ell+\lambda g, \lambda).
$
Equality with $(0,1)$ forces $\lambda=1$ and
$
g \cdot \ell+\ell+g=0$ in $\mathcal{A}.
$
Apply the algebra homomorphism $\Gamma$ (equivalently, take the Gelfand transform to the sides). For every $y \in \mathbb{R}^n$, we get
$$
(\mathscr{H} g)(y)(\mathscr{H} \ell)(y)+(\mathscr{H} \ell)(y)+(\mathscr{H} g)(y)=0 .
$$
Factor the left-hand side to obtain
$
(1+(\mathscr{H} g)(y))(\mathscr{H} \ell)(y)=-(\mathscr{H} g)(y).
$
Due to $1+(\mathscr{H} g)(y) \neq 0$, we may divide and deduce the desired identity. Therefore 
$
(\mathscr{H} \ell)(y)=-\frac{(\mathscr{H} g)(y)}{1+(\mathscr{H} g)(y)}, \forall y \in \mathbb{R}^n.
$ Since $(\ell,1)\in\mathcal A_1$, we deduce that $\ell\in\mathcal A=L_1(\mathbb R^n)$. This completes the proof.
\end{proof}
One thing to note about the Wiener--Lévy invertibility criterion \ref{WienerHalgebra} will be discussed in detail below.
\begin{remark}[Wiener--Lévy, positive form]\label{quantrong}
\textup{Let $\mathscr H: L_1 (\R^n) \to C_0 (\R^n)$	be the transform is defined by \eqref{eq1.3}. Equip $L_1 (\R^n)$ with the product $f \cdot g:=f\underset{\mathscr H}{\ast} g$ is the convolution\eqref{eq2.1}. Suppose $g \in L_1\left(\mathbb{R}^n\right)$ and
	$
	1+(\mathscr{H} g)(y) \neq 0,  \forall y \in \mathbb{R}^n .
	$
	There exists a function $\eta \in L_1\left(\mathbb{R}^n\right)$ such that for every $y \in \mathbb{R}^n$, then
	$
	(\mathscr{H} \eta)(y)=\frac{(\mathscr{H} g)(y)}{1+(\mathscr{H} g)(y)}.
	$	
	By Theorem \ref{WienerHalgebra}, we have the result is immediate from the previous statement that, i.e., produces a function $\ell \in L_1\left(\mathbb{R}^n\right)$ such that	
	$
	(\mathscr{H} \ell)(y)=-\frac{(\mathscr{H} g)(y)}{1+(\mathscr{H} g)(y)}, \forall  y \in \mathbb{R}^n,
	$ together with the elementary observation that the negation map $f \mapsto-f$ is an isometry and a bijective linear automorphism of the Banach algebra $\big(L_1\left(\mathbb{R}^n\right), \underset{\mathscr H}{\ast}\big)$.
	Indeed, let $\ell \in L_1\left(\mathbb{R}^n\right)$ be as above and set $\eta:=-\ell$. Then $\eta \in L_1\left(\mathbb{R}^n\right)$, for each $y$, we have
	$$
	(\mathscr{H} \eta)(y)=(\mathscr{H}(-\ell))(y)=-(\mathscr{H} \ell)(y)=-\left(-\frac{(\mathscr{H} g)(y)}{1+(\mathscr{H} g)(y)}\right)=\frac{(\mathscr{H} g)(y)}{1+(\mathscr{H} g)(y)} .
	$$
	This establishes the asserted existence of the function $\eta$.
It must be emphasized that the two formulations (with the minus sign or without it) are equivalent: one is obtained from the other by replacing the inverse element by its negative. Which of the two is written in a main theorem is merely a matter of convention; it is important only that the algebraic sign is tracked consistently in the algebraic equation one solves in the unitized Banach algebra.}
\end{remark}

Finally, we arrived at the most important theorem of this section. For the reader's convenience, we will present it as a separate subsection.

\subsection{Gelfand’s Spectral Radius Theorem for $\mathscr H$-algebra}
Let $\mathcal{A}=\big(L_1(\mathbb{R}^n),\underset{\mathscr{H}}{\ast},\|.\|_{\a,L_1}\big)$. In the Subsection \ref{sub3.2}, we established the following assertions:
\begin{itemize}
	\item  $\mathcal{A}$ is a commutative Banach algebra under an equivalent norm $\|.\|_{\a,L_1}$ (see Theorem \ref{theorem22}).
	\item Maps
	$\mathscr H: L_1 \to C_0$ is a continuous algebra homomorphism by Lemma \ref{riemann-lebegues}, and $\mathscr H (L_1 (\R^n))$ has density property in $C_0 (\R^n)$ by Proposition \ref{density-H-image}.
\end{itemize}	

 \noindent  For $f \in \mathcal{A}$ define the $k$-th $\underset{\mathscr{H}}{\ast}$-power $f^{\underset{\mathscr{H}}{\ast}}$ inductively by $f^{\underset{\mathscr{H}}{\ast}1}=f$ and $f^{\underset{\mathscr{H}}{\ast}(k+1)}=f \underset{\mathscr{H}}{\ast} f^{\underset{\mathscr{H}}{\ast} k}$. Based on the concept of spectral (see {\bf 3.1}, Chapter VI, p.195 in \cite{Conway90}), since
  $\mathcal{A}$
  be a complex commutative Banach algebra, for any $f\in \mathcal{A}$ we define its spectrum
 $$\sigma_\mathcal{A}(f)=\big\{\lambda \in \C: f- \lambda.1 \ \text{is not invertible in the unitization}\ \mathcal{A}_1 \big\},\ \text{where} \ \mathcal A_1=\mathcal A\oplus\mathbb C, $$ 
  and following the spectral radius (see {\bf 3.7}, Chapter VI, p.197 in \cite{Conway90}) defined by
 $
 r_\mathcal{A}(f):=\sup \left\{|\lambda|: \lambda \in \sigma_\mathcal{A}(f)\right\}.
 $
 \vskip 0.2cm
\noindent {\bf Statement of Theorem.} 
For every $f \in \mathcal{A}$ the following equalities hold:
$$\boxed{
r_\mathcal{A}(f)=\lim _{k \rightarrow \infty}\big\|f^{\underset{\mathscr{H}}{\ast} k}\big\|_{\alpha, L_1(\R^n)}^{\frac{1}{k}}=\sup _{y \in \mathbb{R}^n}|(\mathscr{H} f)(y)|=\|\mathscr{H} f\|_{L_{\infty}\left(\mathbb{R}^n\right)}.
}
$$
\noindent{\bf Proof.} 
Let $\mathfrak{M}(\mathcal{A})$ denote the maximal ideal space (the set of nonzero multiplicative linear functionals on $\mathcal{A}$). For $\varphi \in \mathfrak{M}(\mathcal{A})$ we write $\hat{f}(\varphi)=\varphi(f)$ for the Gelfand transform. By the results proved in Subsection \ref{sub3.2}, every character of $A$ is evaluation at a point: there is a bijection $y \mapsto \Phi_y$ from $\mathbb{R}^n$ onto $\mathfrak{M}(\mathcal{A})$, with $\Phi_y(f)=(\mathscr{H} f)(y)$. Also $\mathscr{H}$ is continuous and $\|\mathscr{H} f\|_{L_\infty} \leq \textup{Const.}\|f\|_{\alpha, L_1}$.
We need to prove the {\bf Statement of Theorem} via three steps (three inequalities)

\noindent{\bf Step 1.} We will show the inequality $\sup_{y \in \mathbb{R}^n}|(\mathscr{H} f)(y)| \leq r_{\mathcal{A}}(f)$.

\noindent  Let $\varphi \in \mathfrak{M}(\mathcal{A})$. Then $\varphi$ is multiplicative, so $\varphi\left(f^k\right)=(\varphi(f))^k$ for each integer $k \geq 1$, where $f^k$ denotes the $k$-th $\underset{\mathscr{H}}{\ast}$-power $f^{\underset{\mathscr{H}}{\ast} k}$. If $|\varphi(f)|>r_\mathcal{A} (f)$, we set $\lambda=\varphi(f)$. By concept of the spectrum, $\lambda \notin \sigma_\mathcal{A}(f)$, hence $(f-\lambda.1)$ is invertible in $\mathcal{A}_1$. But any multiplicative functional $\varphi$ cannot vanish on invertible elements (if $u$ is invertible then $1=\varphi(1)=\varphi(u) \varphi\left(u^{-1}\right)$ so $\varphi(u) \neq 0$ ). In particular $\varphi(f-\lambda. 1)=\varphi(f)-\lambda=$ 0 contradicts invertibility. Thus $$|\varphi(f)| \leq r_\mathcal{A} (f).$$ Taking the supremum over $\varphi \in \mathfrak{M}(\mathcal{A})$, we obtain
$
\sup _{\varphi \in \mathfrak{M}(\mathcal{A})}|\hat{f}(\varphi)| \leq r_\mathcal{A} (f) .
$
Using the identification $\varphi=\Phi_y$ leading to
$$
\sup _{y \in \mathbb{R}^n}|(\mathscr{H} f)(y)| \leq r_{\mathcal{A}}(f) .
$$
\noindent{\bf Step 2.} We will show the inequality $r_{\mathcal{A}}(f) \leq \liminf_{k \to \infty}\big\|f^{\underset{\mathscr{H}}{\ast} k}\big\|_{\alpha, L^1(\R^n)}^{\frac{1}{k}}.
$

\noindent For each integer $k \geq 1$, the spectral mapping (see Theorem 10.28, p.263 in \cite{Rudin1991functional}) gives $\sigma_{\mathcal{A}}(f^k)=\left\{\lambda^k: \lambda \in \sigma_{\mathcal{A}}(f)\right\}.$ Hence $r_{\mathcal{A}}(f^k )=r_{\mathcal{A}}(f)^k.$ On the other hand $\sigma_{\mathcal{A}}\left(f^k\right)$ is contained in the closed disk of radius $\|f^k \|_{\alpha, L_1 (\R^n)}$, because if $|\lambda|>\left\|f^k\right\|_{\alpha, L_1 (\R^n)}$ then the Neumann series argument shows $(f^k-\lambda. 1)$ is invertible (indeed $\big\| \frac{f^k}{\lambda} \big\|<1$ and $\sum_{m \geq 0}\big(\frac{f^k}{\lambda}\big)^m$ is converges). Thus
$
r_{\mathcal{A}}(f)^k=r_{\mathcal{A}}\left(f^k\right) \leq\big\|f^{\underset{\mathscr{H}}{\ast}k}\big\|_{\alpha, L_1(\R^n)}.
$
Taking $k$-th roots and then limit inferior yields
$$
r_{\mathcal{A}}(f) \leq \liminf_{k \to \infty}\big\|f^{\underset{\mathscr{H}}{\ast}k}\big\|_{\alpha, L_1(\R^n)}^{\frac{1}{k}}.
$$

\noindent{\bf Step 3.} We will show the inequality $\limsup _{k \rightarrow \infty}\big\|f^{\underset{\mathscr{H}}{\ast}k}\big\|_{\alpha, L_1(\R^n)}^{\frac{1}{k}} \leq \sup _{y \in \mathbb{R}^n}|(\mathscr{H} f)(y)| .
$

\noindent Because $\mathscr{H}$ is a continuous algebra homomorphism and the Gelfand transform equals $\mathscr{H}$ under the identification $\mathfrak{M}(\mathcal{A}) \simeq \mathbb{R}^n$, we have for each $k$, then
$
\widehat{f^{\underset{\mathscr{H}}{\ast} k}}=\hat{f}^k$ and hence $\|\widehat{f^{\underset{\mathscr{H}}{\ast}k}}\|_{L_\infty(\R^)}=\|\hat{f}\|_{L_\infty(\R^n)}^k=\|\mathscr{H} f\|_{L_\infty(\R^n)}^k .
$
By the general inequality $\|b\| \geq\|\widehat{b}\|_{\infty}$ valid for any $b \in \mathcal{A}$ (because evaluation at each character is a bounded functional of norm $\leq 1$), we then obtain
$$
\|f^{\underset{\mathscr{H}}{\ast}k}\|_{\alpha, L_1(\R^n)} \geq\|\widehat{f^{\underset{\mathscr{H}}{\ast} k}}\|_{L_\infty(\R^n)}=\|\mathscr{H} f\|_{L_\infty(\R^n)}^k.
$$
Taking $k$-th roots, we infer
$
\|f^{\underset{\mathscr{H}}{\ast}k}\|_{\alpha, L_1(\R^n)}^{\frac{1}{k}} \geq\|\mathscr{H} f\|_{L_\infty(\R^n)}.
$
Thus the limit superior of the left-hand side is at least $\|\mathscr{H} f\|_{L_\infty(\R^n)}$, this equivalently to
\begin{equation}\label{limSup}
\limsup _{k \to \infty}\|f^{\underset{\mathscr{H}}{\ast}k}\|_{\alpha, L_1(\R^n)}^{\frac{1}{k}} \geq\|\mathscr{H} f\|_{L_\infty(\R^n)}.
\end{equation}
To get the reversed inequality, observe that if $|\lambda| > \|\mathscr{H} f\|_{L_\infty(\R^n)}$ then $(\mathscr{H} f - \lambda)$ does not vanish on $\mathfrak{M}(\mathcal{A})$. By Gelfand theory (invertibility criterion) \cite{Kaniuth09}, then $(f-\lambda.1)$ is invertible in $\mathcal{A}_1$, hence $|\lambda|> \|\mathscr{H} f\|_{L_\infty(\R^n)}$ is not in $\sigma_\mathcal{A}(f)$. Therefore $\sigma_\mathcal{A}(f)$ is contained in the closed disk $\left\{ |\lambda| \leq \|\mathscr{H} f\|_{L_\infty(\R^n)}\right\}$, so $r_{\mathcal{A}}(f) \leq \|\mathscr{H} f\|_{L_\infty(\R^n)}$. Combining with {\bf Step 1}, which showed $\|\mathscr{H} f\|_{L_\infty(\R^n)} \leq r_{\mathcal{A}}(f)$ yields the equality
\begin{equation}\label{r=H}
r_{\mathcal{A}}(f) = \|\mathscr{H} f\|_{L_\infty(\R^n)}.
\end{equation}
Plugging the equality \eqref{r=H} into {\bf Step 2} gives
\begin{equation}\label{limInf}
\|\mathscr{H} f\|_{L_\infty(\R^n)} \leq \liminf_{k \to \infty}\big\|f^{\underset{\mathscr{H}}{\ast} k}\big\|_{\alpha, L^1(\R^n)}^{\frac{1}{k}}.
\end{equation}
Coupling \eqref{limInf} with the inequality \eqref{limSup} in {\bf Step 3}, then the 
limit inferior and limit superior
 coincide and equal $\|\mathscr{H} f\|_{L_\infty(\R^n)} $. Hence the limit exists and
$$\lim\limits_{k \to \infty}\big\|f^{\underset{\mathscr{H}}{\ast} k}\big\|_{\alpha, L^1(\R^n)}^{\frac{1}{k}}=\|\mathscr{H} f\|_{L_\infty(\R^n)}=r_{\mathcal{A}}(f).$$
This completes the proof of the {\bf Statement of Theorem} for the spectral radius for $f\in \mathcal A$.

\begin{remark}[Great!]
\textup{With the previous identification of the spectrum in the unitization, then
$$(g,1)\  \text{is invertible in}\ \mathcal{A}_1 \Leftrightarrow 1+ (\mathscr H g)(y)\neq 0, \forall y \in \R^n,$$
because  invertibility of $(g,1)$
 is equivalent to the function $1+ \hat{g}$ 
not vanishing on 
$\mathfrak{M}(\mathcal{A})$, i.e., $1+(\mathscr H g)$ 
never vanishes on $\R^n$. This means that the Wiener–Lévy invertibility criterion \ref{WienerHalgebra} can be seen as a corollary of  Gelfand’s spectral theorem for $\mathscr H$-algebras.}

\end{remark}

\section{Young-type theorem and inequality}\label{sec3}
\noindent Similar to Theorem 2.24 in \cite{AdamsFournier2003sobolev} for classical Fourier convolution, we obtain an analogous form for the $\underset{\mathscr{H}}{\ast}$-convolution.
\begin{theorem}[Young's theorem  for $\underset{\mathscr{H}}{\ast}$-convolution]\label{thm3.1}
Let $p,q,r$ are real numbers in open interval $(1,\infty)$  such that 
$1/p +1/q+1/r =2$.
If $f \in L_q(\mathbb{R}^n)$, $g \in L_p(\mathbb{R}^n)$, and $h \in L_r(\mathbb{R}^n)$, then the following inequality holds:
\begin{equation}\label{eq3.1}
\bigg|\int_{\mathbb{R}^n} \big(f \underset{\mathscr{H}}{\ast} g\big)(x) h(x) dx \bigg|
\leq\mathcal{C}\|f\|_{L_q(\mathbb{R}^n)} \, \|g\|_{L_p(\mathbb{R}^n)} \, \|h\|_{L_r(\mathbb{R}^n)},
\end{equation}
where
$
\mathcal{C}	:= \frac{1}{|a|(2\pi)^{n/2}}
\left(\frac{|3a^2-b^2|^q + 3(a^2+b^2)^q}{4}\right)^{1/q}
$. 
\end{theorem}
\begin{proof}
	Let $p_1,q_1,r_1$  be the conjugate exponents of $p,q,r$, respectively, i.e., $1/p + 1/p_1=1$ (and similar for other constants), together with the assumption $1/p +1/q+1/r =2$, the following relations between the conjugate exponents immediately follow
	\begin{equation}\label{eq3.2}
	1/p_1 +1/q_1 +1/r_1 =1
	\quad \text{and}\quad p(1/q_1 +1/r_1) = q(1/p_1 +1/r_1) =r(1/p_1 +1/q_1)=1.
	\end{equation}
Define auxiliary functions $A(x,v):= \big|K_{(a,b)}[f](x,v)\big|^{q/ p_1}\, |h(x)|^{r/p_1}$ and  $B(x,v):= \big|K_{(a,b)}[f](x,v)\big|^{q/r_1}\, |g(v)|^{p/r_1},$ and $C(x,v):= |g(v)|^{p/q_1} \, |h(x)|^{r/q_1}.$
Clearly, $A\in L_{p_1}(\mathbb{R}^{2n})$, $B\in L_{r_1}(\mathbb{R}^{2n})$, and $C\in L_{q_1}(\mathbb{R}^{2n})$.  By \eqref{eq3.2} and the definition \eqref{eq2.1}, we obtain
\begin{align*}
A(x,v) B(x,v) C(x,v) &= \left|K_{(a,b)}[f](x,v)\right|^{q(1/p_1 +1/r_1) }|g(v)|^{p(1/q_1 +1/r_1)} |h(x)|^{r(1/p_1 +1/q_1)}\\
&= \left|K_{(a,b)}[f](x,v)\right| \cdot |g(v)| \cdot |h(x)|.
\end{align*}
Consequently,
$
\big|\int_{\mathbb{R}^n} \big(f \underset{\mathscr{H}}{\ast} g\big)(x)h(x)dx\big|
\;\leq\; \frac{1}{4|a|(2\pi)^{n/2}} 
\int_{\mathbb{R}^{2n}} A(x,v)\, B(x,v)\, C(x,v)\, dx\, dv.
$ Since $\tfrac{1}{p_1}+\tfrac{1}{q_1}+\tfrac{1}{r_1}=1.$ By applying the Hölder inequality, we deduce

\begin{align}
\bigg|\int\limits_{\mathbb{R}^n} \big(f \underset{\mathscr{H}}{\ast} g\big)(x) h(x) dx\bigg|&\leq \frac{1}{4|a|(2\pi)^{n/2}} \bigg\{\int_{\mathbb{R}^{2n}} A^{p_1}(x,v) dx dv\bigg\}^{\frac{1}{p_1}} \bigg\{\int_{\mathbb{R}^{2n}} B^{q_1}(x,v) dx dv\bigg\}^{\frac{1}{q_1}} \bigg\{\int_{\mathbb{R}^{2n}} C^{r_1}(x,v) dx dv\bigg\}^{\frac{1}{r_1}}\nonumber \\
&=\frac{1}{4|a|(2\pi)^{n/2}}\|A\|_{L_{p_1}(\mathbb{R}^{2n})}  \|B\|_{L_{q_1}(\mathbb{R}^{2n})} \|C\|_{L_{r_1}(\mathbb{R}^{2n})}. \label{eq3.4}
\end{align} 
Next we estimate each factor.  
Using \eqref{eq2.3} and Fubini's theorem, we obtain
\begin{align*}
\|A\|_{L_{p_1}(\mathbb{R}^{2n})}^{p_1} &= \int_{\mathbb{R}^{2n}}\big|K_{(a,b)}[f](x,v)\big|^q |h(x)|^r dx dv = \int\limits_{\mathbb{R}^n}|h(x)|^r \bigg(\int_{\mathbb{R}^n}\big|K_{(a,b)}[f](t)\big|^q dt\bigg)dx\\
&=\bigg(\int_{\mathbb{R}^n}\big|K_{(a,b)}[f](t)\big|^q dt\bigg)\bigg(\int_{\mathbb{R}^n}|h(x)|^r dx\bigg)\\
&\leq 4^{q-1}\big(|3a^2-b^2|^q + 3(a^2+b^2)^q\big)\|f\|_{L_q(\mathbb{R}^n)}^q \|h\|_{L_r(\mathbb{R}^n)}^r.
\end{align*}
Therefore,
\begin{equation}\label{eq3.5}
\|A\|_{L_{p_1}(\mathbb{R}^{2n})} 
\leq 4^{\tfrac{q-1}{p_1}}
\big(|3a^2-b^2|^q + 3(a^2+b^2)^q\big)^{1/p_1}
\|f\|_{L_q(\mathbb{R}^n)}^{q/p_1} 
\|h\|_{L_r(\mathbb{R}^n)}^{r/p_1}.
\end{equation}
By an argument analogous to \eqref{eq3.5}, we deduce
\begin{equation}\label{eq3.6}
\|B\|_{L_{r_1}(\mathbb{R}^{2n})} 
\leq 4^{\tfrac{q-1}{r_1}}
\big(|3a^2-b^2|^q + 3(a^2+b^2)^q\big)^{1/r_1}
\|f\|_{L_q(\mathbb{R}^n)}^{q/r_1} 
\|g\|_{L_p(\mathbb{R}^n)}^{p/r_1}.
\end{equation}
	Finally, using Fubini’s theorem again, we have
\begin{align*}
\|C\|_{L_{q_1}(\mathbb{R}^{2n})}^{q_1}&= \int_{\mathbb{R}^{2n}}|g(v)|^{p} |h(x)|^r dv dx =\bigg(\,\int_{\mathbb{R}^n}|g(v)|^p dv\bigg) \bigg(\,\int_{\mathbb{R}^n} |h(x)|^r dr\bigg)= \|g\|_{L_p(\mathbb{R}^n)}^p \|h\|_{L_r(\mathbb{R}^n)}^r,
\end{align*}
which implies
\begin{equation}\label{eq3.7}
\|C\|_{L_{q_1}(\mathbb{R}^{2n})}
= \|g\|_{L_p(\mathbb{R}^n)}^{p/q_1} \,
\|h\|_{L_r(\mathbb{R}^n)}^{r/q_1}.
\end{equation}
Combining \eqref{eq3.5}, \eqref{eq3.6}, and \eqref{eq3.7}, and recalling the relations in \eqref{eq3.2}, we achieve

\begin{align}
&\|A\|_{L_{p_1}(\mathbb{R}^{2n})}
\|B\|_{L_{r_1}(\mathbb{R}^{2n})}
\|C\|_{L_{q_1}(\mathbb{R}^{2n})} \nonumber \\
&\leq 
4^{(q-1)(\tfrac{1}{p_1}+\tfrac{1}{r_1}+\tfrac{1}{q_1})}
\Big(|3a^2-b^2|^q + 3(a^2+b^2)^q\Big)^{\tfrac{1}{p_1}+\tfrac{1}{r_1}+\tfrac{1}{q_1}}\times
\|f\|_{L_q(\mathbb{R}^n)}
\|g\|_{L_p(\mathbb{R}^n)}
\|h\|_{L_r(\mathbb{R}^n)} \nonumber \\
&=
4^{\,1-1/q}
\Big(|3a^2-b^2|^q + 3(a^2+b^2)^q\Big)^{1/q}
\|f\|_{L_q(\mathbb{R}^n)} \|g\|_{L_p(\mathbb{R}^n)} \|h\|_{L_r(\mathbb{R}^n)}. \label{eq3.8}
\end{align}
Substituting \eqref{eq3.8} into \eqref{eq3.4} yields the desired estimate \eqref{eq3.1}, which completes the proof.
\end{proof}
\noindent Comparison with the classical Young inequality. For the classical Fourier convolution, the optimal Young constant is equal to $1$ (for the normalized $L_1$ norm). Our constant $\mathcal{C}$ reduces to $1$ only in certain case, namely: 
 When $a=1, b=0$ (pure cosine transform) and $n=1$ (one-dimension), we get $\mathcal{C} = \frac{1}{\sqrt{2\pi}} \left( \frac{3^q+3}{4} \right)^{1/q}$, which still is  greater than $1$ for $q\geq 1$.
When $a=0, b=1$ (pure sine transform), the formula is not applicable (division by $a$).
This suggests $\mathcal{C}$ is generally larger than the optimal constant. Hence, $\mathcal{C}$ is only an explicit upper-bound on the right-hand side of \eqref{eq3.1}.\\
\noindent In a special case, when the function $h(x)$ is given by \eqref{eq2.1}, then the following Young-type inequality emerges as a direct consequence of Theorem~\ref{thm3.1}. Similar to the Young-type inequality based on the classical Fourier transform \cite{Stein1972Weiss,Neil63}, we will also obtain a consequent
result associated with our convolution for the $\mathscr{H}$-transform.
The argument we employ is inspired by techniques developed in \cite{TuanVKT,PTM25}. 
In particular, our approach follows the strategy outlined in these works, with special emphasis on the choice of an appropriate bounded linear functional, as in \cite{TuanVKT}, which enables an effective application of Riesz's representation theorem.

\begin{corollary}[Young-type inequality]\label{cor3.1}
The $\underset{\mathscr{H}}{\ast}$-convolution \eqref{eq2.1} , is a continuous
bilinear map between suitable $L_r(\mathbb{R}^n)$ space in the sense  that: If $1\leq p,q,r \leq \infty$ satisfy
$
\frac{1}{p}+\frac{1}{q}=1+\frac{1}{r},
$ then $
\underset{\mathscr{H}}{\ast}: L_q(\mathbb{R}^n) \times L_p(\mathbb{R}^n) \to L_r(\mathbb{R}^n)
$  is defined by the impact $(f,g)\mapsto (f \underset{\mathscr{H}}{\ast} g)$ for any functions 
$f\in L_q(\mathbb{R}^n)$ and $g\in L_p(\mathbb{R}^n)$. Moreover, the following estimate always holds:	
	\begin{equation}\label{eq3.9}
	\big\|f \underset{\mathscr{H}}{\ast} g\big\|_{L_r(\mathbb{R}^n)}
	\le \mathcal{C}\|f\|_{L_q(\mathbb{R}^n)}\|g\|_{L_p(\mathbb{R}^n)} \ \text{where}\ \mathcal{C}	:= \frac{1}{|a|(2\pi)^{n/2}}
	\left(\frac{|3a^2-b^2|^q + 3(a^2+b^2)^q}{4}\right)^{1/q}.
	\end{equation}
	
\end{corollary}

\begin{proof}
	We treat three different cases associated with the parameter $r$.
	
	\medskip\noindent\textbf{Case 1.} For $p=q=r=1$.  
	This is the $L_1$-estimate already established in \eqref{eq2.4} of Theorem~\ref{thm2.1}, hence inequality \eqref{eq3.9} holds with the stated constant 	$
	\mathcal{C}$.
	
	\medskip\noindent\textbf{Case 2.} For $1<p,q,r<\infty$.  
	Let $r_1$ denote the Hölder conjugate of $r$, so that $1/r+1/r_1=1$ and $r_1\in(1,\infty)$. From the relation
	$
	\frac{1}{p}+\frac{1}{q}=1+\frac{1}{r}
	$
	we infer
	$
	\frac{1}{p}+\frac{1}{q}+\frac{1}{r_1}=2,
	$  which shows the numbers $p,q,r_1$ satisfies the condition of Theorem~\ref{thm3.1} (with the role of $r$ being replaced by $r_1$). Consequently, for every  $f\in L_q(\mathbb{R}^n)$ and $g\in L_p(\mathbb{R}^n)$, then the linear operator $\mathscr{T}_h:=\int_{\mathbb{R}^n}(f \underset{\mathscr{H}}{\ast} g)(x) h(x) dx$ is bounded continuous in $L_{r_1}(\mathbb{R}^n)$. By applying Riesz's  representation theorem \cite{Stein1972Weiss}, we deduce that $f \underset{\mathscr{H}}{\ast} g$ belongs to $L_{r}(\mathbb{R}^n)$.  To continue the proof,
	we choose $$
	h(x):= \mathrm{sign}\big\{(f \underset{\mathscr{H}}{\ast} g)(x)\big\} \times \big[(f \underset{\mathscr{H}}{\ast} g)(x)\big]^{^{r / r_1}},
	$$ which yields 
		$|h(x)|^{r_1}=|(f \underset{\mathscr{H}}{\ast} g)(x)|^r,$
		and 
		$$\begin{aligned}
\|h\|_{L_{r_1}(\mathbb{R}^n)}=\bigg(\int_{\mathbb{R}^{2n}} |h(x)|^{r_1}dx\bigg)^{1 / r_1}&=\bigg(\int_{\mathbb{R}^{2n}}|(f \underset{\mathscr{H}}{\ast} g)(x)|^r dx\bigg)^{1 / r_1}=\|f \underset{\mathscr{H}}{\ast} g\|_{L_r (\mathbb{R}^n)}^{r / r_1}.
		\end{aligned}$$
And we have an intermediate equality 
$
\int_{\mathbb{R}^n} (f \underset{\mathscr{H}}{\ast} g)(x) h(x) d x=\int_{\mathbb{R}^n}|(f \underset{\mathscr{H}}{\ast} g)(x)|^{\frac{r}{r_1+1}} d x.
$
Since $\frac{1}{r}+\frac{1}{r_1}=1$, then $\frac{r}{r_1+1}=r$. Therefore
\begin{equation}\label{11}
\int_{\mathbb{R}^n}(f \underset{\mathscr{H}}{\ast} g)(x) h(x) d x=\|f \underset{\mathscr{H}}{\ast} g\|_{L_r (\mathbb{R}^n)}^r.
\end{equation}
On the other hand, due to Theorem \ref{thm3.1}, we obtain
\begin{equation}\label{22}
\left|\int_{\mathbb{R}^n}(f \underset{\mathscr{H}}{\ast} g)(x) h(x) d x\right| \leq \mathcal{C}\|f\|_{L_q (\mathbb{R}^n)}\|g\|_{L_p (\mathbb{R}^n)}\|h\|_{L_{r_1} (\mathbb{R}^n)}=\mathcal{C}\|f\|_{L_q (\mathbb{R}^n)}\|g\|_{L_p (\mathbb{R}^n)}\|f\underset{\mathscr{H}}{\ast} g\|_{L_r (\mathbb{R}^n)}^{r / r_1},
\end{equation} where
$
\mathcal{C}	:= \frac{1}{|a|(2\pi)^{n/2}}
\left(\frac{|3a^2-b^2|^q + 3(a^2+b^2)^q}{4}\right)^{1/q}.
$
Coupling \eqref{11} and \eqref{22}, we deduce
$$
\|f \underset{\mathscr{H}}{\ast} g\|_{L_r (\mathbb{R}^n)}^r \leq \mathcal{C}\|f\|_{L_q (\mathbb{R}^n)}\|g\|_{L_p (\mathbb{R}^n)}\|f \underset{\mathscr{H}}{\ast} g\|_{L_r (\mathbb{R}^n)}^{r / r_1}.
$$This equivalent to 
$
\|f \underset{\mathscr{H}}{\ast} g\|_{L_r (\mathbb{R}^n)}^{r-\frac{r}{r_1}} \leq \mathcal{C}\|f\|_{L_q (\mathbb{R}^n)}\|g\|_{L_p (\mathbb{R}^n)}.
$
Since $r-\frac{r}{r_1}=r\big(1-\frac{1}{r_1}\big)=1$, then we direct derive the
inequality \eqref{eq3.9}.

\noindent\textbf{Case 3.} When $r = \infty$ and $p,q\in (1,\infty)$, then $\frac{1}{p}+\frac{1}{q}=1$. We need to show that
\begin{equation}\label{eq3.11a}
\|f \underset{\mathscr{H}}{\ast} g\|_{L_\infty(\mathbb{R}^n)} \leq \mathcal{C}\|f\|_{L_q(\mathbb{R}^n)} \|g\|_{L_p(\mathbb{R}^n)}.
\end{equation}
Indeed, apply Hölder's inequality with pairs exponents $q$ and $p$ and using kernel estimate \eqref{eq2.3}, we obtain
$$\begin{aligned}
\|f \underset{\mathscr{H}}{\ast} g\|_{L_\infty(\mathbb{R}^n)} &\leq \mathrm{ess \sup}_{x\in \mathbb{R}^n} \frac{1}{4|a|(2\pi)^{\frac{n}{2}}} \int_{\mathbb{R}^n} |K_{(a,b)}[f](x,v)|\ |g(v)| dv\\
&\leq \mathrm{ess \sup}_{x\in \mathbb{R}^n} \frac{1}{4|a|(2\pi)^{\frac{n}{2}}}\bigg(\int_{\mathbb{R}^n}|K_{(a,b)}[f](x,v)|^q dv\bigg)^{\frac{1}{q}} \|g\|_{L_p(\mathbb{R}^n)}= \mathcal{C}\|f\|_{L_q(\mathbb{R}^n)} \|g\|_{L_p(\mathbb{R}^n)}
\end{aligned}$$
\end{proof}

\begin{remark}
		\textup{In the one-dimensional case with $a=b=1$, the inequalities \eqref{eq3.1} and \eqref{eq3.9} reduce to estimates (2.5) and (2.14) respectively in  \cite{tuan22}. Hence, this result can be viewed as a generalization of Theorems 2.3 and 2.4 in \cite{tuan22}.}	
\noindent \textup{From Definition \ref{def2.1}, the $\underset{\mathscr{H}}{\ast}$-convolution can be expressed as
$
		(f \underset{\mathscr{H}}{\ast} g)(x) 
		= \frac{1}{4a}\big\{(3a^2-b^2)I_1 + (a^2+b^2)\big(I_2+I_3-I_4\big)\big\},
$ where by setting
		\[
		\begin{aligned}
		I_1 &:= \frac{1}{(2\pi)^{n/2}}\int_{\mathbb R^n} f(x-v)g(v)\,dv, 
		&\ I_2 &:= \frac{1}{(2\pi)^{n/2}}\int_{\mathbb R^n} f(x+v)g(v)\,dv, \\
		I_3 &:= \frac{1}{(2\pi)^{n/2}}\int_{\mathbb R^n} f(-x+v)g(v)\,dv, 
		&\ I_4 &:= \frac{1}{(2\pi)^{n/2}}\int_{\mathbb R^n} f(-x-v)g(v)\,dv.
		\end{aligned}
		\]
\noindent Clearly, $I_1$ coincides with the classical Fourier convolution. Let $\check g$ denote the reflection of $g$, i.e.,
	$\check g(x):=g(-x),\ x\in\mathbb R^n.$
		If $f\in L_p(\mathbb R^n)$ and $g\in L_q(\mathbb R^n)$, then
	$$I_2 = \frac{1}{(2\pi)^{n/2}}\int_{\mathbb R^n} f(x+v)g(v)\,dv
		= \frac{1}{(2\pi)^{n/2}}\int_{\mathbb R^n} f(x-v)\check g(v)\,dv
		= (f\underset{F}{*}\check g)(x),$$
		where $\underset{F}{*}$ denotes the Fourier convolution. Note that $\check{g} \in L_q (\R^n)$ and $\|\check{g}\|_{L_q (\R^n)}=\|g\|_{L_q (\R^n)}$ for any $q\geq 1.$  Hence, by Young’s inequality for the Fourier convolution \cite{AdamsFournier2003sobolev}, we have
$$\begin{aligned}
\|I_2\|_{L_r (\R^n)} &=\|(f\underset{F}{*}\check{g})(x)\|_{L_r (\R^n)} \\&\leq \|f\|_{L_p (\R^n)} \|\check{g}\|_{L_q (\R^n)} = \|f\|_{L_p (\R^n)} \|g\|_{L_q (\R^n)}.
\end{aligned}$$ 
Similarly, $\|I_1 \|_{L_r (\R^n)}=\|(f\underset{F}{*}g)(x)\|_{L_r (\R^n)} \leq \|f\|_{L_p (\R^n)} \|g\|_{L_q (\R^n)}.$ By the same reasoning, we get 
$$\begin{aligned}\|I_3 \|_{L_r (\R^n)}=\|(f\underset{F}{*}\check{g})(-x)\|_{L_r (\R^n)}&=\|(f\underset{F}{*}\check{g})(x)\|_{L_r (\R^n)}\\&\leq \|f\|_{L_p (\R^n)} \|\check{g}\|_{L_q (\R^n)} = \|f\|_{L_p (\R^n)} \|g\|_{L_q (\R^n)},\end{aligned}$$  and $\|I_4 \|_{L_r (\R^n)}=\|(f\underset{F}{*}g)(-x)\|_{L_r (\R^n)}=\|(f\underset{F}{*}g)(x)\|_{L_r (\R^n)}\leq \|f\|_{L_p (\R^n)} \|g\|_{L_q (\R^n)}.$	
Thus, for $p,q,r\ge1$ satisfying
		$
		\frac{1}{p}+\frac{1}{q}=1+\frac{1}{r},
		$
		the triangle inequality yields estimate
		\begin{align}\label{3.11a'}
		\|f\underset{\mathscr H}{*}g\|_{L^r(\mathbb R^n)}
		&\le \frac{|3a^2-b^2|^r+3(a^2+b^2)^r}{4|a|}
		\|f\underset{F}{*}g\|_{L_r(\mathbb R^n)}\le \frac{|3a^2-b^2|^r+3(a^2+b^2)^r}{4|a|}
		\|f\|_{L_p(\mathbb R^n)}\|g\|_{L_q(\mathbb R^n)}.
		\end{align}
It follows that an explicit upper bound in \eqref{eq3.9} is better than that in \eqref{3.11a'}, and hence Corollary \ref{cor3.1} is genuinely non-trivial.}
\end{remark}
\section{Issues related to applicability}\label{sec4}

\subsection{Fredholm integral equation}
In this part, we study a class of integral equations associated with the convolution structure defined in \eqref{eq2.1}. Our objective is to derive sufficient conditions under which the existence and uniqueness
of closed-form solutions can be guaranteed. Following \cite{Sri92}, we consider the Fredholm integral equation of the form
\begin{equation}\label{eq:fredholm}
f(x) + \lambda \int_{\R^n} \Phi_f(x,v)\,g(v)\,dv = h(x),\ x\in\R^n,
\end{equation}
where $\lambda\in\C$, $f$ is an unknown function to be determined, $\Phi_f$ is given kernel and $g,h$ are prescribed functions. We will show the conditions for solvability in $L_1$ of  equation \eqref{eq:fredholm} for case the kernel $
\Phi_f(x,v)=\frac{1}{4a(2\pi)^{n/2}}\,K_{(a,b)}[f](x,v),
$ where $K_{(a,b)}[f]$ is the kernel appearing in \eqref{eq2.2}. Taking $\lambda=1$ and 
$
h(x)=(g\underset{\mathscr H}{\ast} k)(x)
$. By definition of $\underset{\mathscr H}{\ast}$ in 
 \eqref{eq2.1}, then \eqref{eq:fredholm} can be rewritten in convolution form
\begin{equation}\label{eq:fredholm-H}
f(x) + (f\underset{\mathscr H}{\ast} g)(x) \;=\; (g\underset{\mathscr H}{\ast} k)(x),
\ x\in\R^n.
\end{equation}
\begin{theorem}\label{thm:fredholm-solve}
	Let $g,k \in L_1 (\R^n)$ are given functions. For the equation \eqref{eq:fredholm-H} to be solvable in the space $L_1(\R^n)$, a sufficient condition is that $1 + (\mathscr H g)(y)\neq 0$ for every $y\in\R^n.$ Under this condition, the solution $f\in L_1(\R^n)$ exists uniquely and is given almost everywhere on $\R^n$ by the formula $f(x)= (\eta \underset{\mathscr H}{\ast} k)(x)$, where $\eta \in L_1$ is defined via  $$(\mathscr H\eta)(y) = \frac{(\mathscr H g)(y)}{1+(\mathscr H g)(y)}.$$ Moreover, the solution $f$ satisfies the
	estimate $\|f\|_{L_1(\R^n)} \le 
	\frac{|3a^2-b^2|+3(a^2+b^2)}{4|a|(2\pi)^{n/2}}\;\|\eta\|_{L_1(\R^n)}\;\|k\|_{L_1(\R^n)}.$
\end{theorem}

\begin{proof}
	Applying the $\mathscr H$-transform on both sides of	\eqref{eq:fredholm-H}, we obtain 
	$\mathscr H\big[f + (f\underset{\mathscr H}{\ast} g)\big] \;=\; \mathscr H(g\underset{\mathscr H}{\ast} k).$ By utilizing the factorization identity 
	\eqref{eq2.5}
	and the linearity of $\mathscr H$-transform, this leads to the pointwise identity
	\[
	(\mathscr H f)(y)\bigl(1+(\mathscr H g)(y)\bigr) \;=\; (\mathscr H g)(y)\,(\mathscr H k)(y) \ \text{for all}\ y\in\R^n,
	\] which means that \[
	(\mathscr H f)(y) \;=\; \frac{(\mathscr H g)(y)}{1+(\mathscr H g)(y)}\;(\mathscr H k)(y).
	\]
	By the non-vanishing assumption, i.e., $1 + (\mathscr H g)(y)\neq 0, \forall y\in\R^n,$ then the above expression is valid. By Remark \ref{quantrong} and also algebraic structure was developed in Theorem \ref{WienerHalgebra} leading to there exists a function $\eta \in L_1 (\R^n)$ such that $(\mathscr H\eta)(y) = \frac{(\mathscr H g)(y)}{1+(\mathscr H g)(y)}$. This yields $$(\mathscr H f)(y) \;=\; (\mathscr H\eta)(y) (\mathscr H k)(y)=\mathscr H (\eta \underset{\mathscr H}{\ast} k)(y),$$
	which implies, by the injectivity of $\mathscr H$-transform (Lemma \ref{donanh}), then 	
	$f(x)= (\eta \underset{\mathscr H}{\ast} k)(x)$ almost everywhere on $\R^n$. Moreover, since both $\eta$ and $k$ belong to $L_1(\R^n)$  the convolution $\eta \underset{\mathscr H}{\ast} k$ is well-defined and belongs
	to $L_1(\R^n)$; see Theorem~\ref{thm2.1}. The uniqueness of the solution follows from the fact that the convolution structure $\underset{\mathscr H}{\ast}$ is uniquely defined on $L_1$. Therefore, the $L_1$-estimate of the solution $f$ is obtained directly from inequality \eqref{eq2.4}. This completes the proof.
\end{proof}
\noindent We now examine several special cases concerning the boundedness of the solution that arise from the
application of the inequalities established in Section \ref{sec3}. By inequalities \eqref{eq3.9} and \eqref{eq3.11a}, combining the same reasoning with the Young inequality proved in Section \ref{sec3} yields $L_r$-bounds for the solution $f$ when $\eta\in L_q$, $k\in L_p$ and the indices satisfy the usual relation $1/p+1/q=1+1/r$.  In particular, if $p,q\in(1,\infty)$ and $1/p+1/q=1$ (so $r=\infty$), then
\[
\|f\|_{L_\infty(\R^n)} \le \mathcal C\,\|\eta\|_{L_q(\R^n)}\|k\|_{L_p(\R^n)},
\]
where the constant $\mathcal C$ may be taken equal to the explicit upper-bound constant appearing in Corollary~\ref{cor3.1}.

\begin{example}
\textup{Consider in 1-dimension $(n=1)$ with $a=1$ and $b\in \R / \{0\}$. Take $g(x)=k(x)=\sqrt{\frac{\pi}{2}}e^{-|x|} \in L_1 (\R)$.\\ We have $(\mathscr{H}g)(y)= \frac{1}{1+y^2}$ (the value is independent of the parameter $b$ because the sine integral is canceled out due to the evenness of the Gaussian function). Therefore, the sufficient condition in Theorem \ref{thm:fredholm-solve} is $1+(\mathscr H g)(y)$ 
	never vanishes on $\R$ always vaild, and we obtain $(\mathscr H\eta)(y) = \frac{(\mathscr H g)(y)}{1+(\mathscr H g)(y)}=\frac{1}{2+y^2}$. Apply the inverse transform \eqref{eq1.4}, we infer that $\eta(x)=\frac{\sqrt{\pi}}{2} e^{-\sqrt{2}|x|}\in L_1 (\R)$. This yields $(\mathscr H f)(y) \;=\; \mathscr H (\eta \underset{\mathscr H}{\ast} k)(y)=(\mathscr H\eta)(y) (\mathscr H k)(y)=\frac{1}{(2+y^2)(1+y^2)}.$ Thus, the solution 
	$f(x)  =\mathscr{H}^{-1}\left[\frac{1}{1+y^2}-\frac{1}{2+y^2}\right](x) 
	=\sqrt{\frac{\pi}{2}} e^{-|x|}-\frac{\sqrt{\pi}}{2} e^{-\sqrt{2}|x|}$ belongs to $L_1 (\R)$.}
\end{example}

\subsection{A Cauchy problem for the heat equation in one dimension}
We study the initial-value problem
	\begin{align}
	k\frac{\partial^2u(x,t)}{\partial x^2}&=\frac{\partial u(x,t)}{\partial t},\quad (x,t)\in (\mathbb{R}\times\mathbb{R}^+),\label{heat-pde}\\
	u(x,0)&=\varphi(x),\quad x\in\mathbb{R},\label{heat-ic}
	\end{align}
where \(k>0\) is the diffusion coefficient. We assume that, for each \(t>0\),
$
\lim_{|x|\to\infty} u(x,t)=0=\lim_{|x|\to \infty}\frac{\partial u(x,t)}{\partial x}$,
so that the boundary terms arising in integration by parts vanish. The initial datum \(\varphi(x) \) will be taken in \(L_1(\R)\) (or in other Lebesgue classes when required).  In this final part, we will derive the solution to the above problem and represent it using the $\underset{\mathscr{H}}{\ast}$-convolution \eqref{eq2.1}.

\noindent First, recall that in one dimension the \(\mathscr H\)-transform \eqref{eq1.3} is defined by
$$
(\mathscr H f)(y) = \frac{1}{\sqrt{2\pi}}\int_{\R}\bigl(a\cos(xy)+b\sin(xy)\bigr)f(x)\,dx,
\ y\in\R,
$$
with the coefficients \(a,b\) not both zero.
For each fixed \(t>0\), we define \(U(y,t):=(\mathscr H u(x,t))(y)\). Because the \(\mathscr H\)-transform acts only in the \(x\)-variable we may interchange it with the time derivative and obtain
$$
\frac{\partial U(y,t)}{\partial t}=\mathscr H\bigg(\frac{\partial u(x,t)}{\partial t}\bigg)(y).
$$
We now compute \(\mathscr H\big(\frac{\partial^2u(x,t)}{\partial x^2}\big)\). For smooth \(u(x,t)\) that together with $\frac{\partial u(x,t)}{\partial x}$ vanishes at infinity we may integrate by parts twice: using the even/odd symmetry of the cosine and sine kernels one checks that the boundary terms vanish and that
\[
\mathscr H\bigg(\frac{\partial^2u(x,t)}{\partial x^2}\bigg)(y) = -y^2\,\mathscr H\big(u(x,t)\big)(y) = -y^2 U(y,t).
\]
Hence applying \(\mathscr H\)-transform to both-sides of \eqref{heat-pde} and \eqref{heat-ic} yields, for each fixed \(y\in\R\), the ordinary differential equation
\begin{equation}
\Bigg\{
\begin{aligned}
\frac{\partial U(y,t)}{\partial t} &= -ky^2 U(y,t),\quad (y,t)\in (\mathbb{R}\times\mathbb{R}^+),\\
U(y,0)&=(\mathscr{H}\varphi)(y),\quad y\in\mathbb{R}.
\end{aligned}\label{eq4.15}
\end{equation}
The ODE \eqref{eq4.15} is linear and separable; its unique solution is
\begin{equation}\label{U-solution}
U(y,t) = e^{-k t y^2}\,(\mathscr H\varphi)(y),\quad y\in\R,\; t \in \R^+.
\end{equation}
Moreover, by definition \(u(x,t)=\mathscr H^{-1}U(x,t)\). Thus, the solution of the initial-value problem can be written in the transform form
$
u(x,t) = \mathscr H^{-1}\big( e^{-k t y^2}\,\mathscr H\varphi\big).
$
Now, we represent it as a convolution \eqref{eq2.1}. Indeed, since
\(a\neq0\), we may express the Gaussian factor \(e^{-k t y^2}\) itself as an \(\mathscr H\)-transform of an explicit Gaussian in the \(x\)-variable. Define the time-dependent Gaussian
\[
g_t(\xi):=\frac{1}{\sqrt{k t}}\,e^{-\frac{\xi^2}{4k t}},\quad \xi\in\R,\; t>0.
\]
We compute its cosine transform by using the standard Gaussian Fourier integral, and we get
\[
\int_{\R} e^{-a\xi^2}e^{i y\xi}\,d\xi = \sqrt{\frac{\pi}{a}}\, e^{\frac{-y^2}{4a}},\quad \textup{Re}\ a>0,
\]
take \(a=\tfrac{1}{4k t}\) and obtain
\[
\int_{\R} e^{-\xi^2/(4k t)} e^{i y\xi}\,d\xi = 2\sqrt{\pi k t}\; e^{-k t y^2}.
\]
Hence, since the sine integral of an even integrand vanishes, then
\[
\begin{aligned}
F_c[g_t](y)
&:=\frac{1}{\sqrt{2\pi}}\int_\R \cos(\xi y)\,g_t(\xi)\,d\xi
= \frac{1}{\sqrt{2\pi}}\frac{1}{\sqrt{k t}}\cdot 2\sqrt{\pi k t}\; e^{-k t y^2} = \sqrt{2}\; e^{-k t y^2}.
\end{aligned}
\]
Consequently
$
\mathscr H[g_t](y)= a\,F_c[g_t](y) + b\,F_s[g_t](y) = a\sqrt{2}\; e^{-k t y^2},
$
because \(F_s[g_t]=0\) and \(a\neq0\), then 
\begin{equation}\label{Gauss-H}
e^{-k t y^2} \;=\; \frac{1}{a\sqrt{2}}\; \mathscr H[g_t](y)= \mathscr H \bigg(\frac{e^{-\frac{\xi^2}{4k t}}}{\sqrt{k t}}\bigg)(y).
\end{equation}
Combining \eqref{U-solution} and \eqref{Gauss-H}, and using the factorization  
property \eqref{eq2.5},
we deduce that
\[
U(y,t) \;=\; e^{-k t y^2}\,(\mathscr H\varphi)(y)
\;=\; \frac{1}{a\sqrt{2}}\; \mathscr H[g_t](y) (\mathscr H\varphi)(y)
\;=\; \frac{1}{a\sqrt{2}}\; \mathscr H\big(g_t \underset{\mathscr H}{\ast} \varphi\big)(y).
\]
Applying \(\mathscr H^{-1}\)-transform \eqref{eq1.4} when $a\neq 0$, we arrive at the convolution representation
\begin{equation}\label{u-convolution}
u(x,t) \;=\; \frac{1}{a\sqrt{2}} (g_t \underset{\mathscr H}{\ast} \varphi)(x)= \frac{1}{a \sqrt 2}\bigg(\frac{e^{-\frac{\xi^2}{4k t}}}{\sqrt{k t}}\underset{\mathscr H}{\ast}\varphi\bigg)(x).
\end{equation}
By Theorem~\ref{thm2.1}, this formula yields \(u(x,t)\in L_1(\R)\) for every \(t>0\). Note that the representation \eqref{u-convolution} holds whenever \(a\neq0\). In the case, if \(a=0\) (hence \(b\neq0\)), then identity \eqref{Gauss-H} is fails because the sine transform of the even Gaussian \(g_t\) vanishes identically, i.e., \(F_s[g_t]=0\). In this case, the Gaussian factor cannot be expressed as an \(\mathscr H\)-transform of \(g_t\), and one should work with \(u=\mathscr H^{-1}(e^{-k t y^2}\,\mathscr H\varphi)\) directly. 

\noindent We now deduce several norm estimates for \(u(x,t)\) by using  inequalities for \(\underset{\mathscr H}{\ast}\)-convolution established in Section~\ref{sec3}. 
By Theorem~\ref{thm3.1} and Corollary~\ref{cor3.1}, for the one-dimensional case \( (n=1)\), then the upper bound constant is clearly $\mathcal{C}_1=\frac{1}{|a| \sqrt{2\pi}}
\left(\frac{|3a^2-b^2|^q + 3(a^2+b^2)^q}{4}\right)^{1/q}$ when \(p,q,r\in(1,\infty)\). We consider the following cases for the boundedness of solution \eqref{u-convolution}: \\
{(i)} \emph{The \(L_1\)-case.} Take $p=q=r=1$, based on the inequality \eqref{eq2.4} then  $$
\|u\|_{L_1(\mathbb{R})}= \frac{1}{a\sqrt{2}} \|g_t \underset{\mathscr H}{\ast} \varphi\|_{L_1(\mathbb{R})} \leq \frac{1}{a\sqrt{2}}\mathcal{C}_1  \|g_t\|_{L_1(\mathbb{R})} \|\varphi\|_{L_1(\mathbb{R})}.
$$
We have the \(L_1\)-norm of Gaussian integral is independent of \(t\), i.e., $$\|g_t\|_{L_1(\R)} = \int_{\R}\frac{1}{\sqrt{k t}} e^{-\xi^2/(4k t)}\,d\xi
= \frac{1}{\sqrt{k t}}\cdot \sqrt{4\pi k t} = 2\sqrt{\pi}, \ \text{for}\  k>0\ \text{and}\ t>0.$$ Hence, for \(a\neq0\) then  $
\|u\|_{L_1(\R)} \le \frac{|3a^2-b^2|+3(a^2+b^2)}{4|a|^2}\,\|\varphi\|_{L_1(\R)}.
$

\noindent {(ii)} \emph{General \(L_r\)-estimate.}
Let \(p,q,r>1\) satisfy \(\tfrac1p+\tfrac1q = 1+\tfrac1r\).  By \eqref{eq3.9}, we obtain
$
\|u\|_{L_r(\R)}
\le \frac{\mathcal{C}_1}{|a|\sqrt{2}}\;\|g_t\|_{L_p(\R)}\;\|\varphi\|_{L_q(\R)}.
$
The \(L_p\)-norm of the Gaussian \(g_t\) scales in \(t\) as
$\left\|g_t\right\|_{L^p(\mathbb{R})}=2^{\frac{1}{p}}\left(\frac{\pi}{p}\right)^{\frac{1}{2p}}(k t)^{-\frac{p-1}{2 p}}\ \text{for}\  k>0, t>0,$
so the explicit \(t\)-dependence is given by \((k t)^{-(p-1)/(2p)}\).

\noindent {(iii)} \emph{The \(L_\infty\)-case.} If \(p,q>1\) with \(\tfrac1p+\tfrac1q=1\) (so \(r=\infty\)), then
$
\|u\|_{L_\infty(\R)}
\le \frac{\mathcal{C}_1}{|a|\sqrt{2}}\;\|g_t\|_{L_p(\R)}\;\|\varphi\|_{L_q(\R)},$ and the factor \(\|g_t\|_{L_p (\R)}\) is as above.
\vskip 0.4cm

\noindent \textbf{Acknowledgment} \\
\noindent The author expresses their sincere gratitude to the anonymous referees for carefully reading the manuscript and providing insightful and constructive feedback, which has significantly improved both the rigor and presentation of the paper.\\
\noindent \textbf{Disclosure statement}\\
\noindent No potential conflict of interest was reported by the author(s).\\
\noindent \textbf{Ethics declarations}\\
\noindent We confirm that all the research meets ethical guidelines and adheres to the legal requirements of the study country.\\
\noindent \textbf{Funding}\\ \noindent This research received no specific grant from any funding agency.\\
\noindent \textbf{ORCID}\\
\noindent \textsc{Trinh Tuan}  \url{https://orcid.org/0000-0002-0376-0238}

\end{document}